\documentclass[preprint,12pt]{elsarticle}

\usepackage{lineno,hyperref}
\usepackage{graphicx}
\usepackage{{ntheorem}}
\usepackage[T1]{fontenc}
\usepackage{appendix}
\usepackage{amssymb}
\usepackage[dvipsnames]{xcolor}
\usepackage{amsfonts}
\usepackage{gensymb}
\usepackage{multicol}
\usepackage{anyfontsize}
\usepackage[utf8]{inputenc}
\usepackage{calligra}
\usepackage{xcolor}
\usepackage{pdfpages}
\usepackage{nicefrac}
\usepackage{subfig}
\usepackage{amsmath}
\usepackage{natbib}
\usepackage{esint}
\usepackage{algorithm}
\usepackage{algorithmic}
\usepackage{comment}
\usepackage{mathtools}

\modulolinenumbers[5]

\journal{Journal of \LaTeX\ Templates}

\definecolor{GioGreen}{RGB}{0, 140, 80}
\definecolor{manc}{RGB}{255, 100, 0}

\definecolor{AdditionColor}{rgb}{0.0, 0.22, 0.62}

\newcommand{\modify}[1]{\textcolor{AdditionColor}{#1}}

\newcommand{\ds}{\displaystyle}
\definecolor{Ema}{RGB}{250, 142, 0}

\newcommand{\ha}{\frac{1}{2}}

\newcommand{\tF}{\tilde{F}}

\newcounter{dctr}[section]
\numberwithin{equation}{section}
\numberwithin{figure}{section}

\newtheorem{remark}{Remark}[section]

\begin{document}
	
	\begin{frontmatter}
		
		\title{\textbf{\large 
		A semi-implicit finite volume method \\
		for the Exner model of sediment transport}}
		
		\author{S. Avgerinos, M. Castro, E. Macca G. Russo\\
			University of Catania, Universidad de M\'alaga.}
			

    	\begin{abstract}
            The aim of this work is to construct efficient finite volume schemes for the numerical study of sediment transport in shallow water, in the framework of the Exner model \cite{BoscarinoRusso,CastroNieto}.
            In most cases, the velocity related to the sediment is much lower that the fluid velocity, which, in turn, may be much lower that the free-surface wave speed. Explicit methods that resolve all waves require small time steps due to the CFL stability restriction because of fast surface waves. Furthermore, if Rusanov flux is adopted, slow sediment waves may be affected by the large numerical diffusion.  The objective of the present work is to drastically improve the  efficiency in the computation of the evolution of the sediment by treating water waves implicitly, thus allowing much larger time steps than the one required by fully explicit schemes. The goal is reached by suitably semi-implicit schemes obtained by the use of implicit-explicit Runge-Kutta methods. 
    	\end{abstract}
	\end{frontmatter}

\section{Introduction} 
According to \cite{SimpsonCastelltort}, there exist two perspectives in the study
of the interaction between fluid flow and the dynamics of
sediment transport. The first aspect includes the influence
of bed change on the flow dynamics as is studied, for example, in \cite{DutykhDias,Fernandez,tsunamiMalaga,RzadkiewiczMariotti}, where tsunami waves generated by the
submarine sediment movement or topography motion are considered. The second one includes the morphological bed change and sedimentary deposits due to the flow dynamics. The dam-break over erodible bed, for instance, causes morphological changes of topography as shown in \cite{Audusse-Berthon,Cordier-Le,Audusse-Bouchut,GunawanLhebrard,CastroNieto,VanRijn,Castro-FernandezNieto2,VanRijn_2}. The focus of this work is the second perspective. 

Sediments are particles that can be transported by rivers or sea over large distances. Their transport, deposition or erosion, controls the form of rivers. Depending on the water discharge, the geology, the slope of the area and the supply of sediment, rivers exhibit braid, meander or straight patterns. 
The main applications pertain to the prediction of pattern formation in river beds (such as dunes) and the evolution of river morphology. 
In an industrial context, the accumulation of sediment at the bottom of a hydroelectric dam or of a harbor are relevant problems. Its study is particularly useful to prevent the filling of hydroelectric reservoirs, but also to preserve water intakes in rivers.

The hydrodynamic model described by the Saint-Venant system \cite{SaintVenant_2,SaintVenant_1} of the shallow water equation is widely used to study the flow of fluid in rivers, coastal areas, etc. The morphodynamic model described by the Exner equation is used to model the bed load sediment transport, dam breaks, floods, a particular class of submarine landslides. The Exner model includes  a conservation law related to the evolution of the bottom topography due to the action of the fluid. 
It  requires to express the solid flux as a function of the hydrodynamical variables, water depth $h$,  and velocity $u$,  and
some other physical parameters, such as the mean grain diameter, or the porosity of the sediment layer.  In our case we adopt the Grass equation \cite{Grass} to complete the Saint-Venant-Exner system, but other approximations could be consider without increasing the complexity of the numerical scheme presented here.

There is a vast literature on the Saint-Venant–Exner system. In 2011, Cordier et al. \cite{Cordier-Le} compared a splitting method and a coupled method in solving the shallow water system coupled with sediment transport models, in which the splitting method may produce instability related to a bad estimation of the wave speeds of the complete system. Gunawan and Lhébrard \cite{GunawanLhebrard} proposed a hydrostatic relaxation scheme for the shallow water-Exner equations, which can be seen as the hydrostatic reconstruction of relaxation solvers. Staggered schemes for the Exner-shallow water equations was proposed by Gunawan et al. \cite{GunawanEymard}. Liu et al.\ \cite{LiuBeljadid} proposed a coupled method for water flow, sediment transport and bed erosion on triangular meshes. Audusse et al.\ \cite{AudusseChalons} proposed a robust splitting method for the Saint-Venant–Exner equation based on an approximate Riemann solver. Murillo et al.\ \cite{MURILLO20108704} developed a first order explicit reconstruction adopting a Roe-type scheme. A well-balanced central weighted essentially nonoscillatory scheme for the sediment transport model was presented in Qian et al.\ \cite{QianLi}.

When the sediment speed is much slower than the water save speeds, the system becomes stiff, and the treatment by an explicit method
may introduce a strong restriction on the time step, which makes accurate prediction of the sediment evolution very expensive. 
In one space dimension, the Saint-Venant–Exner system modeled by a $3\times 3$ system of partial differential equations in space and time. A detail analysis about the hyperbolicity of the system is presented in \cite{Cordier-Le}. In many interesting applications, when the Froude number, $F_r=|u|/gh$ is small ($F_r\ll 1$) and the interaction between the fluid and the sediment layer is weak, we have that the two wave speeds related to the hydrodynamical component satisfy that $\lambda_1<0$ and $\lambda_3>0$, and $|\lambda_2|\ll \min(|\lambda_1|,|\lambda_3|),$ i.e. the wave speed of the sediment is much smaller than those of the hydrodynamical component. An explicit method implies a strong stability restriction due to the velocity of the free-surface wave. 
If the evolution of the sediment does not depend on the detailed behaviour of the fast surface waves, then by filtering them out by means of an implicit treatment one should be able to use much larger time step, still obtaining a detailed description of the sediment motion. 

The objective of the present paper is to drastically improve the  efficiency in the computation of the evolution of the sediment by treating water waves implicitly, thus allowing much larger time steps than the one required by explicit schemes.

Other authors used implicit methods for the Exner model.
Rosetti et al. \cite{RosettiBonaventura} have proposed a semi-implicit semi-lagrangian scheme by adopting a $\theta-$method with an appropriate $\theta$ increasing stability hence efficiency. Successively, Bonaventura et al. \cite{BonaventuraMultilayer} presented several semi-implicit, multi-layer methods for the shallow water model and the Exner model. They, in addition to the $\theta$-method, adopt an IMEX-ARK2 method with 3 stages and showed that the introduction of an implicit part makes the method faster, in the computational sense, than the explicit one. Recently, Garres-Díaz et al. (2022) proposed a semi-implicit $\theta-$method approach for sediment transport models \cite{GarresFernandez} by which, choosing $\theta>\ha$, an increase in both efficiency and stability was obtained by Casulli in\cite{CasulliCattani}. 

At variance with what has been previously done, here we propose a semi-implicit method based on IMEX, and at the same time show that the long time evolution of the sediment depends weakly on the detailed behaviour of the surface waves, which opens to the possibility of filtering out fast surface waves, and focus on the sediment flow.
To this purpose, we compare the numerical solution of the Exner model with the evolution of the simple wave associated to the sediment, which may be considered an exact solution of the Exter system, before shock formation. 

The paper is structured in the following way: in Section 2 the one-dimensional model equations are introduced and discussed; in Section 3 the details of the semi-implicit numerical method are presented while the forth section is devoted to the derivation of the scalar equation that defines the simple wave. In Section 5 we present several tests to numerically assess the accuracy of the method and the capability to capture the sediment transport without resolving the detailed evolution of the fast surface waves.  
Section 6-8 are devoted to the two dimensional extension and testing of the method. 
Finally, in section 9 we draw some conclusions.

\section{1D Exner Model}
	Let us consider the one-dimensional hyperbolic shallow water system with bathymetry 
	\begin{equation}
    \label{Shallow_water}    
    \begin{cases}
        h_t + q_x =0\\
     \displaystyle   q_t + \left(\frac{q^2}{h} + \frac{g}{2}h^2 \right)_x = -ghb_x,
    \end{cases}
    \end{equation}
    \begin{figure}[!ht]
	\centering	
	\includegraphics[scale = 0.5]{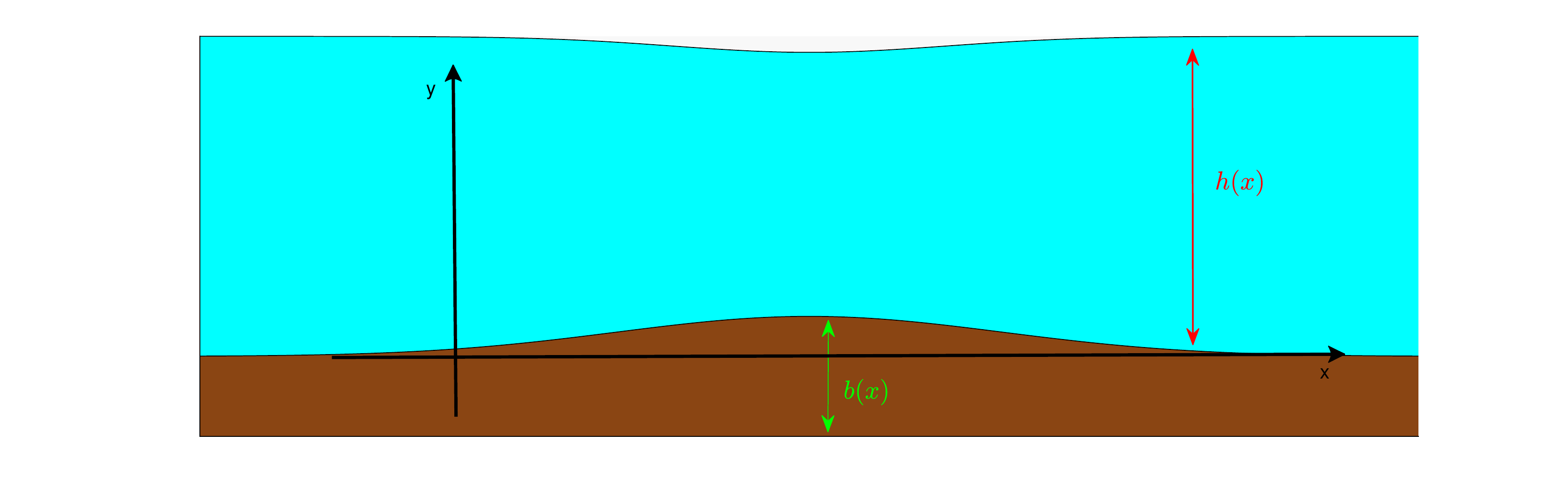}
	\vspace{-0.4 cm}
	\caption{Shallow water equations: water-flow $h(x)$ and bottom topography $b(x).$}
	\label{Fig_Ex_1}
\end{figure}
where $x$ denotes the space coordinate along the axis of the channel and $t$ is time; $q(x,t)$  represents the water flux per unit width (discharge) and $h(x,t)$ the water thickness; $g$ the acceleration due to gravity;  $b(x)$ denotes the bottom topography; furthermore, the following relation holds $q(x,t) =h(x,t)u(x,t),$ with $u$ the depth average horizontal velocity as shown in Figure \ref{Fig_Ex_1}.
    
The system of equations used in this work is obtained by coupling shallow water equation $\eqref{Shallow_water}$ and the sediment equation: 
\begin{equation}
        \label{sediment_equation}
        (z_b)_t + (q_b)_x = 0
\end{equation}
where $z_b(x,t)$ represents the height of sediment layer and  $q_b(h,q)(x,t),$ denotes the solid transport discharge, in our case computed by the Grass model \cite{Grass,QianLi,CastroNieto} 
\begin{equation}
	    \label{q_b}
	    q_b = \xi A_g u|u|^{m-1}
\end{equation}
with $m\in[1,4],$ $A_g\in]0,1[$ and $\xi = {1}/{(1-\rho_0)}$ where $\rho_0$ is the porosity of the sediment layer.
Throughout this paper we shall assume that the porosity is constant.

In this way, the Exner 1D system is given by:
\begin{equation}
	    \label{Ex_sis}
	    \begin{cases}
	        h_t + q_x = 0,\\
	        \displaystyle q_t + \left(\frac{q^2}{h} + \ha gh^2\right)_x = -gh(b + z_b)_x, \\
	        (z_b)_t + (q_b)_x = 0.  
        \end{cases}
\end{equation}
Note that, if $S$ is defined as $S(x,t) = b(x) + z_b(x,t),$ we have $\dfrac{\partial S}{\partial t} = \dfrac{\partial z_b}{\partial t},$ so system \eqref{Ex_sis} could be rewritten as 
\begin{equation}
	    \label{Ex_sis_S}
	    \begin{cases}
	        h_t + q_x = 0,\\
	        \displaystyle q_t + \left(\frac{q^2}{h} + \ha  gh^2\right)_x = -ghS_x, \\
	        S_t + (q_b)_x = 0.  
        \end{cases}
\end{equation}
System $\eqref{Ex_sis}$ can be written as a hyperbolic system  with a non-conservative term 
\begin{equation}
	    \label{non_cons_term}
	    \frac{\partial U}{\partial t} + \frac{\partial F(U)}{\partial x} = B(U)\frac{\partial U}{\partial x},
\end{equation}
where 
\[
	    U = \begin{bmatrix} h \\ q \\ S     \end{bmatrix}, \quad 
	    F = \begin{bmatrix}\ds  q \\ \ds \frac{q^2}{h} + \ha gh^2 \\ \ds q_b     \end{bmatrix},
	    \quad B(U) = \begin{bmatrix} 0 & 0 & 0 \\ 0 & 0 & -gh \\ 0 & 0 & 0     \end{bmatrix},
\]
and $q_b$ is given by eq.  \eqref{q_b}.\\
Given $J = \nabla_UF$ and $A(U) = J(U) - B(U),$ system \eqref{non_cons_term} can be rewritten as 
\begin{equation}
	    \label{non_cons_term_1}
	    \frac{\partial U}{\partial t} + A(U)\frac{\partial U}{\partial x} = 0,
\end{equation}
where,
\[
    A(U) = \begin{bmatrix} 0 & 1 & 0 \\gh - u^2 & 2u & gh \\ \alpha & \beta & 0     \end{bmatrix},
\]
in which
$\alpha = \dfrac{\partial q_b}{\partial h}$ and $\beta = \dfrac{\partial q_b}{\partial q}$.
Assuming $u>0$ in the whole domain one has $\beta = m \xi A_g u^{m-1}/h$ and $\alpha = -u\beta$.

	
This system is strictly hyperbolic if and only if the characteristic polynomial: 
\[
    p_{\lambda}(\lambda) = -\lambda((u-\lambda)^2 -gh) + gh\beta(\lambda -u)
\]
has three distinct real roots $\lambda_1<\lambda_2<\lambda_3.$

As $A_g \rightarrow0,$ $\beta$ vanishes, and the three eigenvalues become $\lambda_1 = u-c$, 
$\lambda_2 = 0$ and $\lambda_3 = u+c$, with 
$c=\sqrt{gh}$. 

We are mainly interested in regimes for which the local Froude number $F_r = |u|/c$ is relatively small, say $F_r<1/2$, and the quantity $\beta$ is much smaller than 1.

For sufficiently small values of $A_g$, such that $\beta\ll 1$, performing an asymptotic expansion of the eigenvalues, 
one obtains 
\begin{align}
\label{lamd_1}
    \lambda_1 &= u - \sqrt{gh} - \beta\frac{\sqrt{gh}}{2(1 - F_r)} + O(\beta^2) \\ \label{lamd_2}
    \lambda_2 &= \frac{\beta g h u}{gh - u^2} + O(\beta^2) = \beta u /(1-F_r^2) + O(\beta^2) \\ \label{lamd_3}
    \lambda_3 &= u + \sqrt{gh} + \beta\frac{\sqrt{gh}}{2(1 + F_r)} + O(\beta^2)   
\end{align}

The parameter $\beta$, therefore, represents a first approximation of the ratio between the sediment wave speed and the fluid speed:
\begin{equation}
    \label{eq:lambda2}
    \frac{\lambda_2}{u} = \frac{\beta}{1-F_r^2} + O(\beta^2)
\end{equation}
	



In the next section we shall derive the semi-implicit numerical method that we adopt in the paper, in which, under the assumption of small Froude number, we treat implicitly the gravity waves of the hydrodynamical component and explicitly the sediment wave.

\begin{figure}[!ht]
	\centering	
	\includegraphics[scale = 0.5]{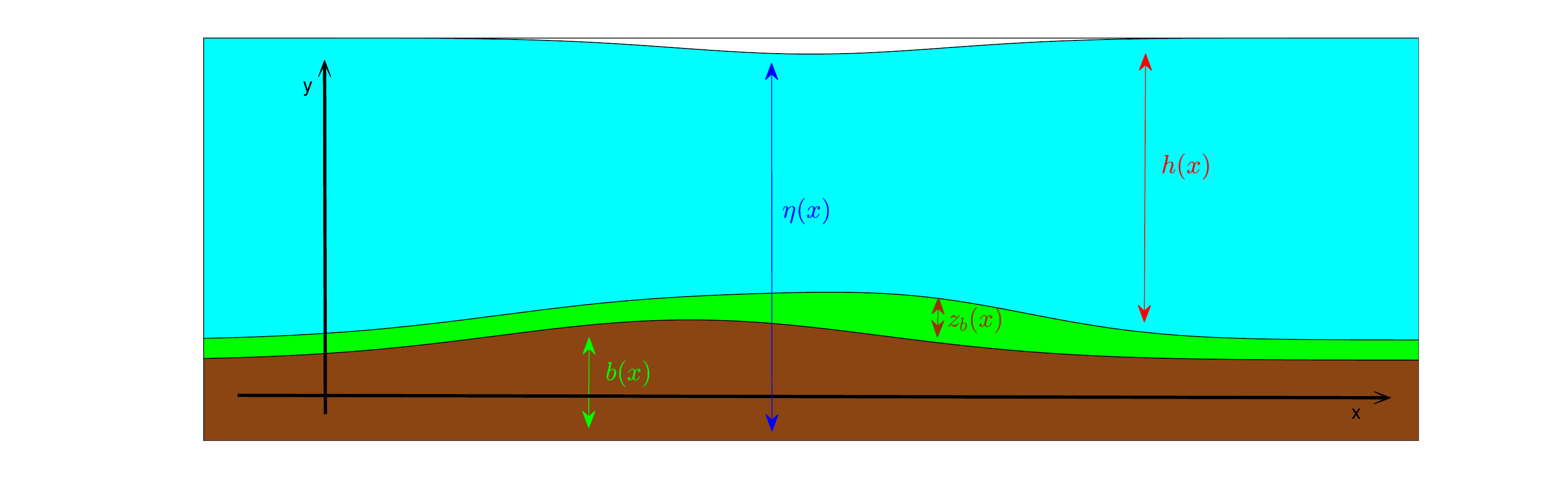}
	\vspace{-0.4 cm}
	\caption{1D Exner model: water surface $\eta(x);$ water-flow $h(x);$ sediment layer $z_b(x)$ and bottom topography $b(x).$}
	\label{Fig_Ex_2}
\end{figure}

It is convenient to  rewrite the 1D-Exner model \eqref{Ex_sis} in terms of  $\eta(x,t) = h(x,t) + b(x) + z_b(x,t) $, which  represents the elevation of the undisturbed water surface, in place of the water thickness $h$ (see Figure \ref{Fig_Ex_2}), because for lake at rest $\eta$ is constant, and for stationary solutions with low Froude numbers $\eta$ is almost constant.  
Thus, system \eqref{Ex_sis_S}, will be written as follows:
	\begin{equation}
	    \label{Ex_sis_eta}
	    \begin{cases}
	        \eta_t + (q+q_b)_x = 0\\
	        q_t + (qu)_x + gh(\eta)_x = 0 \\
	        (z_b)_t + (q_b)_x = 0  
        \end{cases}
	\end{equation}

	\section{Semi implicit scheme} \label{1D semi implicit}
	In this section, we will focus on the introduction of a scheme derived from an implicit treatment of the gravity surface water waves, while the slow wave corresponding to the sediment evolution is treated explicitly. In particular, we will illustrate a first and second order semi-implicit schemes. 
	
	Let us consider a partition of the interval $[a,b]$ in cell defined by $I_i = [x_{i-\ha},x_{i+\ha}],$ with $i = 1,\ldots,N$ For the sake of simplicity, from now on we suppose that all the cells have the same length $\Delta x$ and $x_{i} = a + (i-1/2)\Delta x$ denotes  center of cell $I_i$. Let $\Delta t$ be the time step such that $t^n = n\Delta t.$\footnote{In practice it is better to assign $\Delta t$  dynamically at each time step by imposing a CFL condition. The choice of constant time step here is adopted in order to simplify the notation in the description of the method.}
	
	Definitely, we denote by $U_i^n$ an approximation on the mean value of $U$ over cell $I_i$ at time $t=t^n,$
	$$U^n_i\cong \frac{1}{\Delta x}\int_{x_{i-\ha}}^{x_{i+\ha}}U(x,t^n)dx. $$
	
	\subsection{First order scheme}
	Following the idea proposed in \cite{Boscarino-Filbet} in which IMEX  Runge-Kutta methods are used for systems in which the stiffness is not necessarily of additive or partitioned type, we  write system \eqref{Ex_sis_eta} as a large system of ODE's, in which we adopt suitable discrete operators for the approximation of space derivatives. The key point in \cite{Boscarino-Filbet} is to identify which specific term has to be treated implicitly and which can be treated explicitly. 
    
    Following \cite{Boscarino-Filbet}, we write system \eqref{Ex_sis_eta} in the form
    \begin{equation}
        \label{ode_bosca_semi}
        U' = H(U_E,U_I).
    \end{equation}
    where $U = [\eta, q, z_b]^T$ and $H(U_E,U_I)$ is given by
    \begin{equation}
        \label{ode_form_semi}
        H(U_E,U_I) = 
        \begin{bmatrix}
            -(q_I+(q_b)_E)_x \\
            -((qu)_E)_x - gh_E(\eta_I)_x \\
            -((q_b)_E)_x
        \end{bmatrix}
    \end{equation}
where the subscript $E$ and $I$ denote which term has to be treated explicitly and which implicitly. 

The semi-implicit scheme can be written in the form
    \begin{equation}
        \label{ode_bosca}
        U' = \tilde{H}(U_E,U_I).
    \end{equation}
with
    \begin{equation}
        \label{ode_form_part}
        \tilde{H}(U_E,U_I) = 
        \begin{bmatrix}
           -\hat{D}_x((q_B)_E)  & -D_x(q_I) \\ -\hat{D}_x((qu)_E) &  -gh_ED_x(\eta_I) \\ -\hat{D}((q_b)_E) &
        \end{bmatrix}.
    \end{equation}
	The semi-discrete in time first order semi-implicit scheme can be written as:
	\begin{equation}
	    \begin{cases}
	        \label{First order}
	        \eta^{n+1} = \eta^n - \Delta t \hat{D}_x(q_b^n) - \Delta t D_x(q^{n+1}),\\
	        q^{n+1} = q^n - \Delta t \hat{D}_x(q^nu^n) - \Delta tgh^nD_x(\eta^{n+1}), \\
	        z_b^{n+1} = z_b^n - \Delta t \hat{D}_x(q_b^n),
	    \end{cases}
	\end{equation}
	where the differential operators $D_x$ and $\hat{D}_x$ applied to a given flux function $\tF(U)$ are respectively defined as:
	\begin{itemize}
	    \item $D_x(\tF_i) = \frac{\tF_{i+\ha} - \tF_{i-\ha}}{\Delta x},$ in which $\tF_{i\pm\ha}$ is suitably defined on cell edges; 
	    \item $\hat{D}_x(\tF_i) = \frac{\tF_{i+\ha} - \tF_{i-\ha}}{\Delta x},$ where  $\tF_{i+\ha} = \ha\Bigl( \tF(U_{i+\ha}^{-}) + \tF(U_{i+\ha}^{+}) - \alpha_{i+\ha}\big(U_{i+\ha}^{+} - U_{i+\ha}^{-}\big)\Bigr)$ is the Rusanov flux and $\alpha_{i+\ha}$ is related to the eigenvalues of the explicit sub system. In our case  $\alpha \approx |u|\ll\max(|\lambda_1|,|\lambda_3|).$
	\end{itemize}
 
The values $U^\pm_{i+\ha}$ at the intercells are computed component-wise throughout conservative linear reconstruction with slope limiter obtained by the generalized MinMod function: 
$$ v_{i\pm\ha}^{\mp} = \Bar{v_i} \pm v_{i}'\frac{\Delta x}{2}$$
with $$v_i' = \frac{1}{\Delta x}{\rm MM}\Bigl(\theta(\Bar{v}_i - \Bar{v}_{i-1}), \frac{\Bar{v}_{i+1} - \Bar{v}_{i-1}}{2},\theta(\Bar{v}_{i+1} - \Bar{v}_{i})\Bigr),$$
where $\theta\in [1,2],$ in our case $\theta = 1.9,$ and 
$$ {\rm MM}(a,b,c) = \begin{cases}
    sign(a)\min(|a|,|b|,|c|) \quad {\rm if} \; a,b,c \; {\rm have \, the\, same\, sign} \\ 0 \;\quad\quad \quad \quad \quad \quad \quad \quad \quad \quad {\rm otherwise}.
\end{cases} $$


	For the sake of simplicity, denoting by $\eta^*$ and $q^*$ the explicit part of first and second equation in system \eqref{First order}, it can be rewrite, to distinguish explicit part from implicit one, as:
	\begin{equation}
	    \begin{cases}
	        \label{First order star}
	         q^{*} = q^n - \Delta t \hat{D}_x(q^nu^n);\\
	         \eta^{*} = \eta^n - \Delta t \hat{D}_x(q_b^n) - \Delta t \hat{D}_x(q^{*});\\
	        \eta^{n+1} = \eta^* + g\Delta t^2 D_x(h^nD_x(\eta^{n+1}));\\
	        q^{n+1} = q^* - \Delta tgh^nD_x(\eta^{n+1}); \\
	        z_b^{n+1} = z_b^n - \Delta t \hat{D}_x(q_b^n); \\ 
            h^{n+1} = \eta^{n+1} - z_b^{n+1} - b.
	    \end{cases}
	\end{equation}
	The procedure to solve system \eqref{First order} and consequently \eqref{First order star} is:
	\begin{enumerate}
	    \item compute $q^* = q^n - \Delta t\hat{D}_x(q^nu^n)$ as $$ q^*_i = q_i^n - \frac{\Delta t}{\Delta x}\Bigl(F_{i+\ha,q}^n - F_{i-\ha,q}^n\Bigr), $$ 
	    where $F_{i\pm\ha,q}$ are the Rusanov flux, previously defined, related to the $q;$
	    \item compute $\eta^* = \eta^n - \Delta t\hat{D}_x(q_b^n)-\Delta t\hat{D}_x(q^*)$ as $$ \eta^*_i = \eta^n_i - \frac{\Delta t}{\Delta x}\Bigl( F^n_{{i+\ha,\eta}} - F^n_{{i-\ha,\eta}}\Bigr) - \frac{\Delta t}{2\Delta x}\Bigl( q^*_{i+1} - q^*_{i-1}\Bigr),$$ in which $F_{i\pm\ha,\eta}$ are again computed with the Rusanov flux;
	    \item fixed $k = g\Bigl(\dfrac{\Delta t}{\Delta x}\Bigr)^2,$ solve implicitly $\eta^{n+1} = \eta^* + g\Delta t^2D_x(h^nD_x(\eta^{n+1}))$ in the following way $$\eta_i^{n+1}\bigg(1+k(h_{i+\ha}^n + h_{i-\ha}^n)\bigg) - \eta^{n+1}_{i+1}kh_{i+\ha}^n - \eta^{n+1}_{i-1}kh_{i-\ha}^n = \eta_i^* \quad {\rm{for\; all}} \quad i = 1,\ldots,N,$$
        where $h_{i\pm \ha}^n=\ha \left(h_{i\pm 1}^n+h_i^n \right).$
        
	    This is an invertible tridiagonal linear system which can be solved to detect $\eta^{n+1} = [\eta^{n+1}_1,\ldots,\eta^{n+1}_{N}];$ 
	    \item compute $q^{n+1} = q^* - \Delta tgh^nD_x(\eta^{n+1})$ $$ q_i^{n+1} = q_i^* - \frac{g\Delta t}{\Delta x}h_i^n\Bigl(\eta_{i+\ha}^{n+1}-\eta_{i-\ha}^{n+1}\Bigr),$$ where $\eta_{i\pm\ha}^{n+1} = \ha\Bigl(\eta_{i\pm1}^{n+1} + \eta_{i}^{n+1}\Bigr); $ 
	    \item compute $z_{b}^{n+1} = z_b^n - \Delta t\hat{D}_x(q_b^n)$ as $$ z_{b_i}^{n+1} = z_{b_i}^n - \frac{\Delta t}{\Delta x}\Bigl(F_{i+\ha,z_b}^n - F_{i-\ha,z_b}^n\Bigr), $$ where $F_{i\pm\ha,z_b}$ are computed with the Rusanov flux, in general $F_{i\pm\ha,z_b}^n\neq F_{i\pm\ha,\eta}^n;$
	    \item compute $h^{n+1}_i = \eta^{n+1}_i - b_i - z_{b_i}^{n+1}.$
	\end{enumerate}

\subsection{Second order scheme}
    As have been done for the first order case and following \cite{Boscarino-Filbet}, we write system \eqref{Ex_sis_eta} in the ODE form \eqref{ode_bosca_semi}-\eqref{ode_form_semi}

After space discretization, the semi-implicit scheme can be written in the form
    \begin{equation}
        \label{ode_bosca_or2}
        U' = \tilde{H}(U_E,U_I).
    \end{equation}
with
    \begin{equation}
        \label{ode_form_part_or2}
        \tilde{H}(U_E,U_I) = 
        \begin{bmatrix}
           -\hat{D}_x((q_B)_E)  & -D_x(q_I) \\ -\hat{D}_x((qu)_E) &  -gh_ED_x(\eta_I) \\ -\hat{D}((q_b)_E) &
        \end{bmatrix}.
    \end{equation}
    With this in mind, we apply an IMEX scheme to system \eqref{ode_form_part_or2}.
    The general procedure to update the numerical solution from time $t_n$ to $t_{n+1}$ using an $s$-stage Runge-Kutta IMEX method is the following:
    \begin{itemize}
        \item Stage values: For $i=1,\ldots,s$ compute
        \begin{align*}
            U^{(i)}_E & = U^n + \Delta t\sum_{j=1}^{i-1}a_{i,j}^EH\left(U_E^{(j)},U_I^{(j)}\right)\\
            U^{(i)}_I & = U^n + \Delta t
            \left(\sum_{j=1}^{i-1}a_{i,j}^I H\left(U_E^{(j)},U_I^{(j)}\right) + a_{i,i}^I H\left(U_E^{(i)},U_I^{(i)}\right)\right).
        \end{align*}
        \item Numerical solution:\\
        $U^{n+1} = U_I^{(s)}.$
    \end{itemize}

    \begin{remark}
    Observe that if a system of the form \eqref{ode_bosca_or2} is autonomous (i.e.\ the right hand side does not explicitly depend on time), and the $s$ stage double Butcher tableau has identical $b$ coefficients, then the evolution requires only $s$ evaluation of function $H$. Furthermore, if the last row of matrix $A$ is equal to the weights $b$ (i.e.\ the implicit tableau defines a {\em stiffly accurate} scheme), then the numerical solution coincides with the last stage value of the implicit scheme \cite{Boscarino-Filbet}.
    \end{remark}

Here we consider the IMEX scheme defined by the following  double Butcher tableau \cite{Boscarino-Filbet}: 
\begin{equation}
\label{tableau}
    \begin{array}{c|cc}
             & 0 &  \\
            c & c & 0\\ \hline
            & 1-\gamma & \gamma
        \end{array} 
        \hspace{3 cm}
        \begin{array}{c|cc}
            \gamma & \gamma &  \\
            1 & 1-\gamma & \gamma\\ \hline
            & 1-\gamma & \gamma
        \end{array}
\end{equation}
where $\gamma = 1 - \frac{1}{\sqrt{2}}$ and $c = \frac{1}{2\gamma}.$
In our case, applying the scheme defined by \eqref{tableau} we have:
    \begin{enumerate}
        \item $U_E^{(1)} = U^n;$
        \item $U_I^{(1)} = U^n + \Delta t\gamma H(U_E^{(1)},U_I^{(1)});$
        \item $U_E^{(2)} =  U^n + \Delta tc H(U_E^{(1)},U_I^{(1)});$
        \item $U_I^{(2)} = U^n + \Delta t(1-\gamma) H(U_E^{(1)},U_I^{(1)}) + \Delta t\gamma H(U_E^{(2)},U_I^{(2)});$
        \item $U^{n+1} = U_I^{(2)}.$
    \end{enumerate}
    \begin{remark}
Let observe that $U_E^{(2)},U_I^{(2)}$ and $U_I^{(1)}$ have a common term, thus, step 3 and 4 may be rewritten as:
    \begin{align*}
        U_E^{(2)} &= (1-\frac{c}{\gamma})U^n + \frac{c}{\gamma}U_I^{(1)};  \\
        U_I^{(2)} &= (1-\frac{1-\gamma}{\gamma})U^n + \frac{1-\gamma}{\gamma}U_I^{(1)} +\Delta t\gamma H(U_E^{(2)},U_I^{(2)}).
    \end{align*} 
    \end{remark}

 \paragraph{Stability condition}
 For an explicit scheme, the CFL restriction is determined by the maximum spectral radius of the matrix $A(U)$ that defines the hyperbolic system \eqref{non_cons_term_1}:
 \begin{equation}
    \textrm{CFL} = \frac{\lambda^n_{\rm max} \Delta t_n}{\Delta x} \leq C_{\textrm{ex}}
    \label{CFL_ex}
\end{equation}
where $\lambda^n_{\rm max} = \max_i\rho(A(U^n_j))$, and $C_{\textrm{ex}}$ is a constant close to one. Here $\rho(A)$ denotes the spectral radius of matrix $A$.

For our semi-implicit scheme \eqref{First order star} we empirically find the following stability condition: 
 $$
    \textrm{MCFL} = \frac{u^n_{\rm max} \Delta t_n}{\Delta x} \leq C_{\textrm{im}}
$$
where $u^n_{\rm max} = \max_j{|u^n_j|}$ and $C_{\textrm{im}} \approx 0.85$.

This condition is much less restrictive than condition \eqref{CFL_ex}, since the condition for the classical CFL becomes
\[ 
    \textrm{CFL}\leq \frac{\lambda^n_{\rm max}}{u^n_{\rm max}} C_{\rm im}
\]
and $\lambda^n_{\rm max}/u^n_{\rm max}\gg 1$ for small Froude number. We expect this estimate to be accurate when $\beta\ll 1$, i.e. when the eigenvalues corresponding to the fast waves are close to those of the standard shallow water model.

\section{Scalar Equation for 1D Exner Model} \label{Scal_sect}
Assuming that $A_g \ll 1$, that is for weak coupling,  the motion of the sediment takes place on a much longer time scale than surface waves. For such a reason, surface waves move over a bathymetry given by the bottom and the sediment, which is almost constant in time. We can therefore imagine that to detect the slow motion of the sediment, a reasonable approximation consists in monitoring the sediment motion on a sequence of quasi-stationary states. This is obtained by setting to zero the time derivative in the first two equations of the Exner model. Our starting point is therefore the following: we neglect the time derivative of $\eta$ and $q$ in the first two equations, because these conditions correspond to stationary flow when the sediment does not move.

	\begin{equation}
	    \label{Start_point}
	    \begin{cases}
	        (q+q_b)_x = 0\\
	        (qu)_x + gh(h + z_b + b)_x = 0\\
	        (z_b)_t + (q_b)_x = 0
	    \end{cases}
	\end{equation}
	We shall use the first to equations to express all quantities $\eta$, $q$, and $z_b$ in terms of $u$.
	From the first equation of (\ref{Start_point}) we get 
	\begin{equation}
	    \label{q+q_b=Q}
	    q + q_b = Q
	\end{equation} 
	hence, assuming $u>0,$ 
	\begin{equation}
	    \label{h}
	    h = \frac{Q}{u} - A_g u^{m_g-1}
	\end{equation}
	where, in this section for the sake of simplicity, we include the coefficient $\xi$ in the parameter $A_g.$, and assume we have a net flux, i.e.\ $Q>0$.
	
	Dividing the second equation by $h$, which we assume non zero (without water flux there would be no sediment flux, so $q$ and therefore $h$ are non zero), we obtain: 
	\begin{equation}
	    \label{2_mod}
	    \frac{u}{q}(qu)_x + g(h+z_b+b)_x=0.
	\end{equation}
    Let us define $G(u)$ a function such that 
    \begin{equation}
        \label{G_x}
        \frac{\partial G}{\partial x} = \frac{u}{q}(qu)_x;
    \end{equation}    
    $\dfrac{\partial G}{\partial x} = \dfrac{dG}{d u}u_x $ hence $\dfrac{u}{q}(qu)_x = G'u_x$ and therefore
    \begin{equation}
        G'(u)u_x = \frac{u}{q}(q'u + q)u_x \quad \Rightarrow \quad G'(u) = \frac{q'}{q}u^2 + u.
    \end{equation}
    As a consequence of the first equation of system (\ref{Start_point}) $q' = -q_b' = -m_gA_gu^{m_g-1},$ then $G'$ takes the form 
    \begin{equation}
        \label{G'}
        G'(u) = \frac{Q-(m_g+1)A_gu^{m_g}}{Q - A_gu^{m_g}}u. 
    \end{equation}
    From equation (\ref{2_mod}) we obtain $G+g(h+z_b +b)=C,$ where $C$ is a constant, consequently $z_b = \dfrac{(C-G)}{g} - h -b.$ Furthermore, from the third equation of (\ref{Start_point}) we have $z_b'u_t + q_b'u_x =0$ which implies
    \begin{equation}
        \label{z_b'}
        z_b' = -\frac{G'}{g} + \frac{Q+(m_g-1)A_gu^{m_g}}{u^2}.
    \end{equation}
    Finally, linking all the results obtained, we find the non-linear scalar equation
    \begin{equation}
        \label{Scal_eq}
        u_t + \lambda(u)u_x=0
    \end{equation}
    where 
    \begin{equation}
        \label{lambda}
        \lambda(u) = \frac{m_gA_gu^{m_g-1}}{ \frac{Q + (m_g-1)A_gu^{m_g}}{u^2} - \frac{G'(u)}{g}}.
    \end{equation}
    
Equations (\ref{Scal_eq},\ref{lambda}), together with equations (\ref{q+q_b=Q},\ref{h},\ref{z_b'}), provide a
solution of the sediment transport in the quasi-static approximation. Such solution will loose validity after shock formation.
\begin{remark}
    Notice that this approximation is close, but not equivalent the simple wave corresponding to the second eigenvalue $\lambda_2$.
\end{remark}

\begin{remark}
The value of $\lambda$ \eqref{lambda} must be compared with the approximate value of the eigenvalue $\lambda_2.$ Indeed, with reference to Eq.~\eqref{eq:lambda2},
$$ 
    \lambda(u) = \frac{m_gA_gu^{m_g-1}}{ \frac{Q + (m_g-1)A_gu^{m_g}}{u^2} - \frac{G'(u)}{g}} = \frac{m_gA_gu^{m_g-1}u/h}{ \frac{Q + (m_g-1)A_gu^{m_g}}{uh} - \frac{G'(u)}{g}u/h} \approx \frac{\beta u}{1-F_r^2} \approx \lambda_2
$$ 
where, since $Q = q+q_b,$
$$\frac{Q + (m_g-1)A_gu^{m_g}}{uh} = \frac{q + q_b + (m_g-1)q_b}{q} = \frac{q + m_gq_b}{q}\approx 1 $$ and 
$$\frac{G'(u)}{g}u/h = \frac{\frac{Q-(m_g+1)A_gu^{m_g}}{Q - A_gu^{m_g}}u}{g}u/h = \frac{Q-(m_g+1)A_gu^{m_g}}{Q - A_gu^{m_g}}\frac{u^2}{gh} = \frac{q - m_g q_b}{q} F_r^2 \approx F_r^2. $$ 
\end{remark}

\begin{figure}
    \centering
    \includegraphics[width=0.48\textwidth]{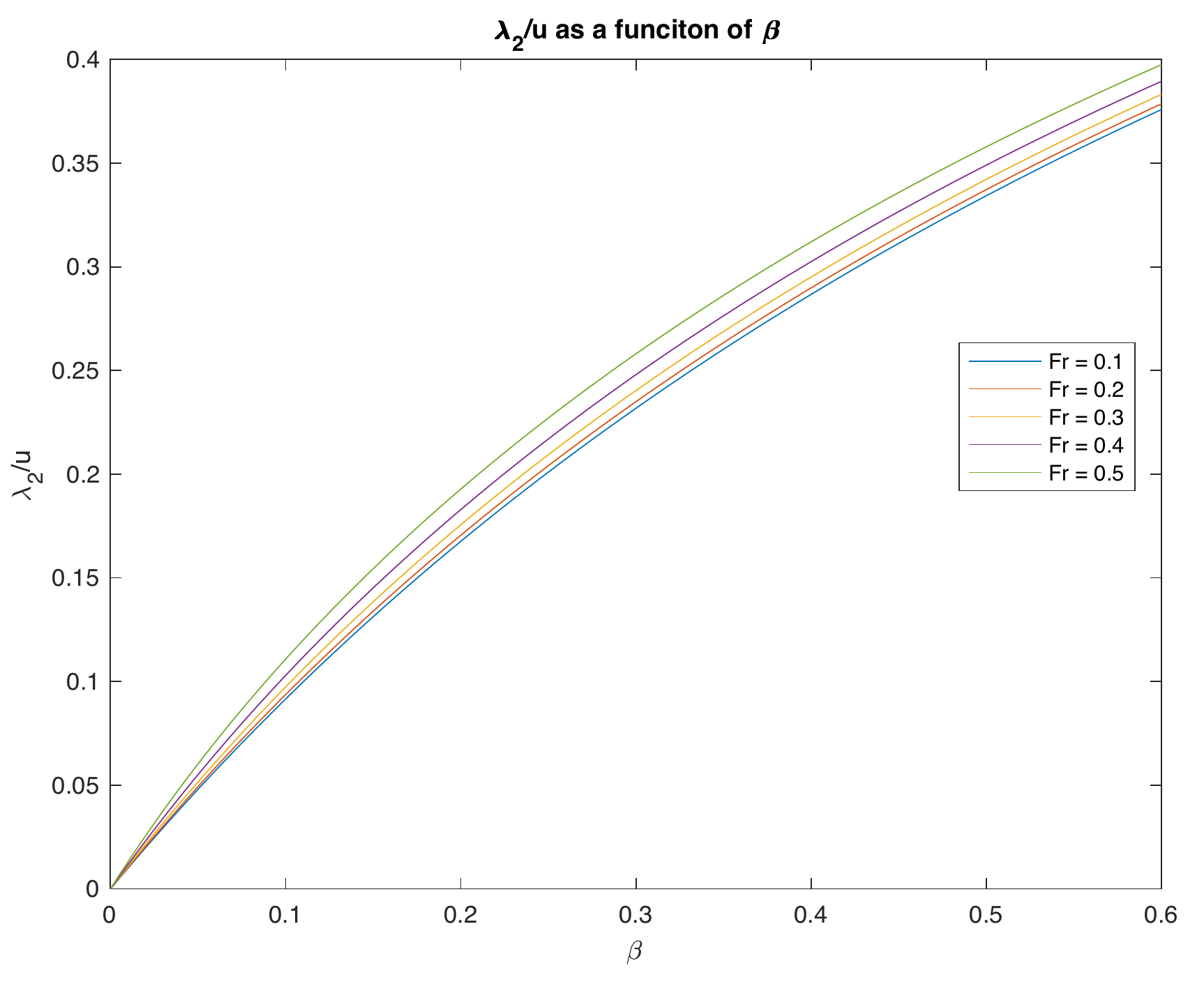}
    \includegraphics[width=0.50\textwidth]{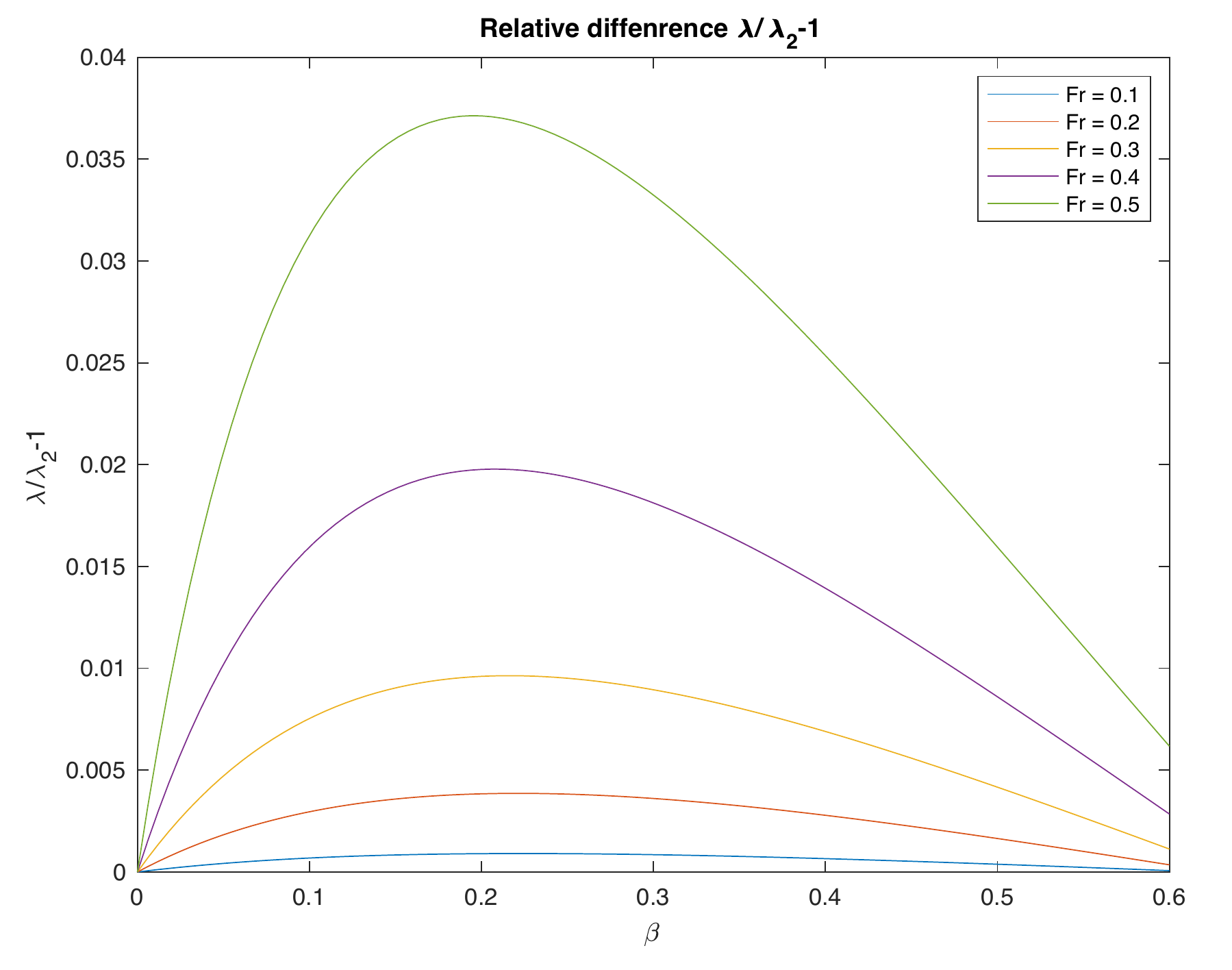}
    \caption{$\lambda_2/u$ as a function of $\beta$ for various Froude numbers (left panel), and relative difference between $\lambda(u)$ and $\lambda_2$ (right panel)} 
    \label{fig:compare}
\end{figure}

After some algebra and knowing that $\beta u = m_gq_b/q$, $\lambda(u)$ can be written as:
\begin{equation}
    \lambda(u) = \frac{\beta u }{(1+\beta) - F_r^2(1-\beta)} 
    = \frac{\beta u}{1-F_r^2+\beta (1+F_r^2)}
\end{equation}
This expression is very close to the one of the second eigenvalue $\lambda_2$. One can show that 
\begin{align*}
    \lambda & = \frac{\beta u}{1-F_r^2}\left(1-\beta\frac{1+F_r^2}{1-F_r^2}\right) + O(\beta^3)\\
    \lambda_2 & = \frac{\beta u}{1-F_r^2}\left(1-\beta\frac{1+F_r^2}{(1-F_r^2)^2}\right) + O(\beta^3)
\end{align*}
therefore the two propagation speeds agree for small values of $\beta$ and Froude number.
The actual difference between $\lambda(u)$ and $\lambda_2$ is actually much smaller than the difference between the two approximate expressions (see Figure \ref{fig:compare}).

\subsection{Second order numerical scheme for the approximated scalar equation}
    In order to solve numerically equation (\ref{Scal_eq}) we adopt the Lax-Wendroff scheme applied to equation in form (\ref{Scal_eq}) \cite{Ricthmeyer_Morton,Carlos}. 
   
    
    The Lax-Wendroff scheme applied to 
    \eqref{Scal_eq} becomes 
    \begin{align}
        \label{L-W}
        u_i^{n+1} = &u_i^n - \Delta t\lambda(u_i^n)\Bigl(\frac{u_{i+1}^n - u_{i-1}^n }{2\Delta x}\Bigr) + &\nonumber \\ & + \frac{\Delta t^2}{2} \lambda(u_i^n)\Bigg[ \lambda'(u_i^n)\Bigl(\frac{u_{i+1}^n - u_{i-1}^n }{2\Delta x}\Bigr) + \frac{\lambda_{i+\ha}^n\Bigl(u_{i+1}^n - u_{i}^n\Bigr) - \lambda_{i-\ha}^n\Bigl(u_{i}^n - u_{i-1}^n\Bigr)}{\Delta x^2} \Bigg].
    \end{align}
    We plan to integrate the equation only in conditions in which the solution remains smooth, and for this reason we shall not use any limiter.

    \section{Numerical experiments}
	When large time steps are used, we should check whether we are able to correctly follow the sediment evolution even if the details on the fast water waves are lost. For this reason, we check the ability of the scheme, presented in previous section, to compute the bathymetry evolution. The main purpose is increase the CFL-condition as much as possible in order to reduce the computational cost using the IMEX strategy described before.
	
    \subsection{1D order accuracy}\label{Test_1} 
    In this section we will compare the solutions obtained by a second order explicit scheme obtained with a third order CWENO reconstruction in space \cite{Russocweno1} and a second order Runge-Kutta scheme \cite{RK2}, first and second order semi-implicit schemes applied to system \eqref{Ex_sis_eta} and first and second order explicit scheme applied to the scalar equation \eqref{Scal_eq}. With this purpose in mind, we have to consider different CFL conditions, and consequently different time steps, one for each scheme that we adopt. For this reason, we define CFL$_{scal}$ the CFL condition used for the scalar equation scheme (\ref{Scal_eq}); CFL$_{expl}$ the CFL condition adopted for the explicit scheme applied to full system  \eqref{Ex_sis_eta}; and CFL$_{IMEX}$ the CFL condition adopted for the semi-implicit scheme applied to full system  \eqref{Ex_sis_eta}. In particular we set CFL$_{scal}=0.9$ for the explicit scalar schemes; CFL$_{expl}=0.4$ for the second order explicit scheme applied to \eqref{Ex_sis_eta}; and  for the semi-implicit methods a larger CFL condition could be used, however, since the term $qu$ \eqref{First order} is treated explicitly, the semi-implicit CFL condition could not be arbitrary larger and a material CFL condition must be satisfied. In our case CFL$_{IMEX} = 15$ is adopted. These CFL conditions come out from the natural stability conditions of the different methods. 
    
    \begin{figure}[!ht]
    	\centering	
    	\hspace{-1.2cm}
    	\includegraphics[scale = 0.48]{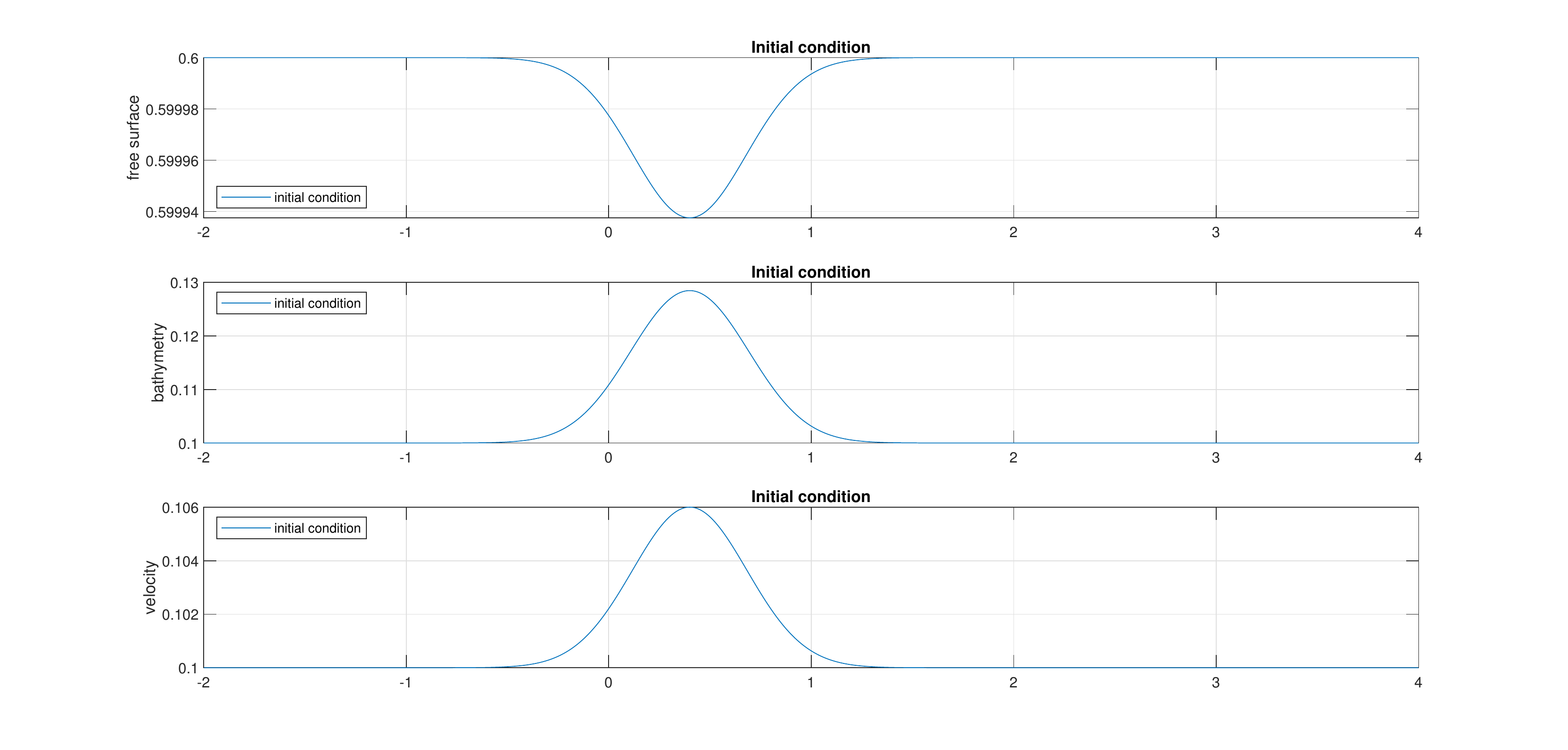}
    	\vspace{-0.8 cm}
    	\caption{Test \ref{Test_1}: (1D order accuracy). Initial condition of free surface (top), sediment layer (center) and velocity (down) for the Exner model on the interval $[-2,4]$ using a $200-$mesh points.  }
    	\label{Test_1_1}
    \end{figure}
    \begin{figure}[!ht]
    	\centering	
    	\hspace{-1.2cm}
    	\includegraphics[scale = 0.16]{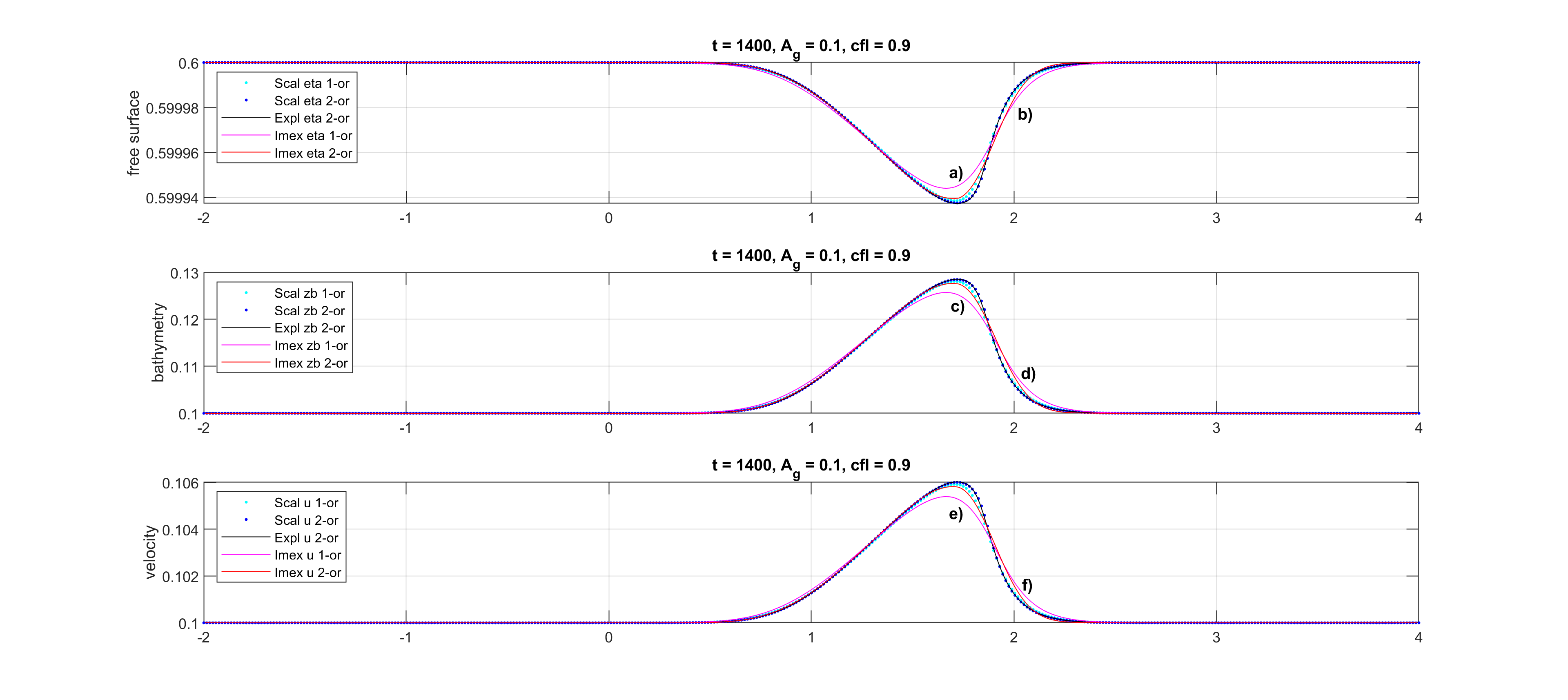}
    	\vspace{-0.8 cm}
    	\caption{Test \ref{Test_1}: (1D order of accuracy). Numerical solutions of free surface (up), sediment layer (center) and velocity (down) for the Exner model on the interval $[-2,4]$ using a $200-$mesh points at time $t=1400$ with, respectively, CFL$_{scal}=0.9,$ CFL$_{expl}=0.4$ and CFL$_{IMEX} = 15.$  }
    	\label{Test_1_2}
    \end{figure}
    \begin{figure}[!ht]
    	\centering	
    	\hspace{-1.2cm}
    	\includegraphics[scale = 0.48]{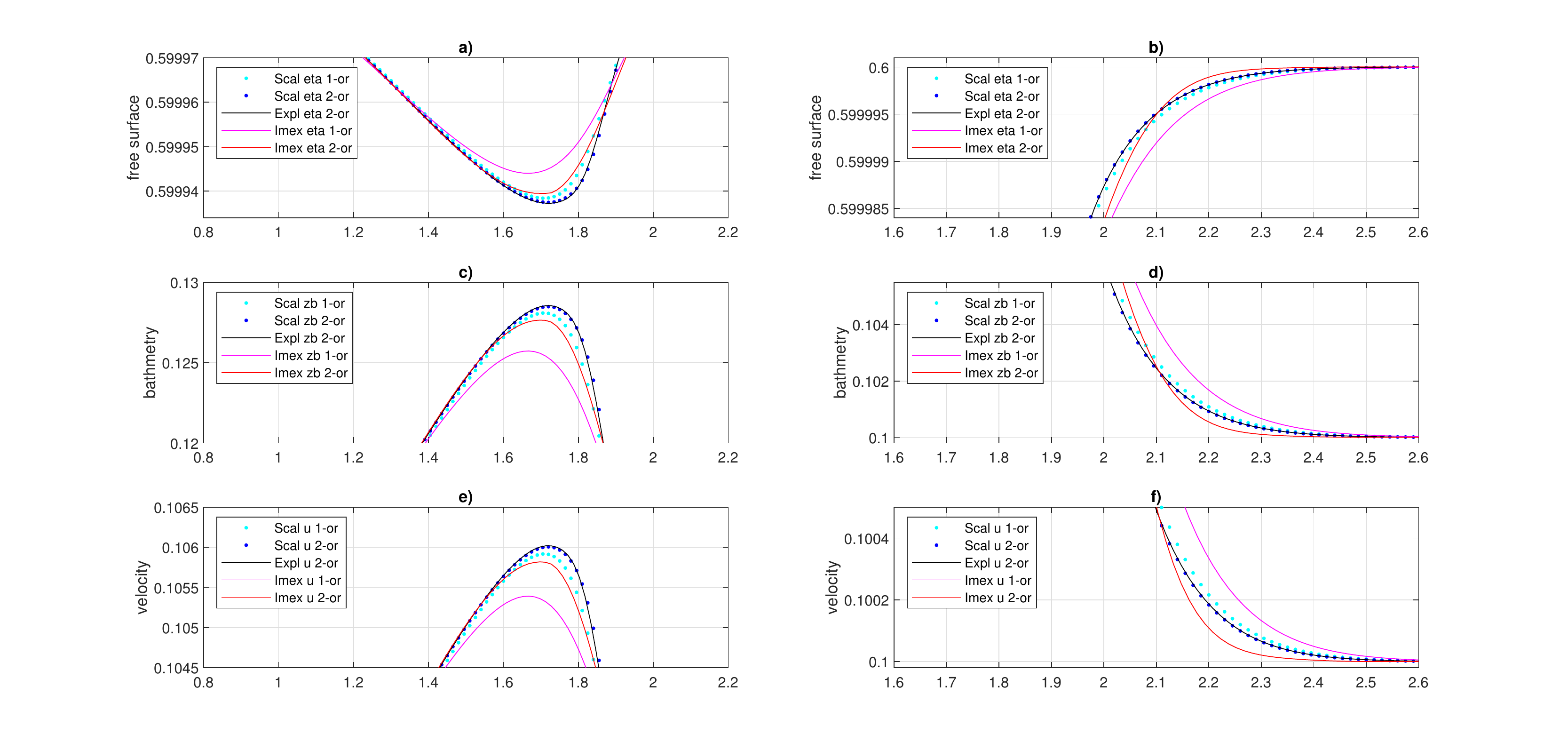}
    	\vspace{-0.8 cm}
    	\caption{Test \ref{Test_1}: (1D order accuracy). Zoom of critical parts for numerical solutions of free surface (up), sediment layer (center) and velocity (down) for the Exner model at time $t=1400$ with, respectively, CFL$_{scal}=0.9,$ CFL$_{expl}=0.4$ and CFL$_{IMEX} = 15.$  }
    	\label{Test_1_3}
    \end{figure}
  
    The common settings of this experiment are: $[a,b]=[-2,4]$ the interval; $A_g = 0.1;$, $\xi = \frac{1}{1-\rho_0},$ with $\rho_0=0.2;$ $m=3;$ $t_{end}= 1400;$ and, since in Section \ref{Scal_sect} all the variables are written depending on the velocity $u,$ initial conditions are so set: $b(x) \equiv 0,$ $h_0(x) = 0.5,$
    
    \begin{equation}
        \label{u0} 
        u_0(x) = 0.1 + 0.006e^{-\frac{(x-0.4)^2}{0.4^2}}
    \end{equation}
    and $z_b(a) = 0.1.$
    
    The constant $Q$ is obtained through $Q = q_0(a) + q_b(a);$ while $C$ is computed as $C = G(u_0(a)) + g(h_0(a) + z_b(a) + b(a))$ where $G$ is a solution of (\ref{G'}) and $g$ is the gravitational constant $g=9.81.$ Free boundary conditions are imposed in both boundaries. 
    
    \begin{table}[htbp]
        \begin{center}
            \begin{tabular}{|c|c|c|c|c|c|c|c|c|c|c|}
                \hline      $z_b$                &\multicolumn{2}{c|}{IMEX or1}&\multicolumn{2}{c|}{IMEX or2}&               \multicolumn{2}{c|}{Scal or1} &  \multicolumn{2}{c|}{Scal or2} & \multicolumn{2}{c|}{Expl or2}  \\
                N   &  Ord & Error      & Ord & Error        & Ord  & Error     & Ord& Error    & Ord& Error    \\ 
                \hline
                200      & -    &  4.41E-4     & -     & 5.34E-4      & -      & 1.31E-4   & -    &  3.09E-5 & -    &  8.00E-5 \\   
                \hline  
                400      & 0.78 &  2.58E-4     & 1.64  & 1.71E-4      & 0.91   & 6.99E-5   & 1.98 &  7.85E-6 & 1.94 &  2.01E-5 \\   
                \hline  
                800      & 0.84 &  1.44E-4     & 2.31  & 3.44E-5      & 0.96   & 3.58E-5   & 2.02 &  1.94E-6 & 2.02 &  5.16E-6 \\   
                \hline  
                1600     & 0.90 &  7.70E-5     & 2.29  & 1.83E-6      & 0.98   & 1.83E-6   & 2.01 &  4.83E-7 & 2.00 &  1.29E-6 \\   
                \hline
            \end{tabular}
            \vspace{2mm}
            \caption{Test \ref{Test_1}: (1D order \modify{of} accuracy) Errors in $L^1-$norm and convergence rates related to the sediment $z_b$ for scalar, explicit and semi-implicit scheme at time $t=1400$ with, respectively, CFL$_{scal}=0.9,$ CFL$_{expl}=0.4$ and CFL$_{IMEX} = 15.$}
            \label{table_1}
            \end{center}
    \end{table}
    
    Figure \ref{Test_1_1} shows the initial condition of the height sediment layer $z_b$ (center), the free surface $\eta$ (up) and velocity $u$ (down). The initial condition of thickness marks out by equation (\ref{h}) while the initial condition of sediment layer comes out from $z_b = (C-G(u_0))/g -h_0(x) -b(x).$ Figure \ref{Test_1_2} exhibits the time evolution of free surface $\eta$ (up);  sediment layer $z_b$ (center); and velocity $u$ (down) at the final time $t = 1400.$ The zoom of critical parts are shown in Figure \ref{Test_1_3}. Table \ref{table_1} proves that all the methods are able to keep the expected order refining the mesh-grid. 
The final time is chosen before a shock forms, since after shock formation the scalar equation does not approximate the solution of the system any more. As expected, there is very good agreement between the solutions of the scalar equation and those obtained with the full system. This is due to the good quality of the quasi-stationary approximation adopted to derive the scalar equation, which consequently provides a good approximation of a simple wave of the system.

    In this setting of conditions it is impossible to understand if the discrepancy between the two solutions (scalar and system) is  due to the numerical approach or the errors introduced for the derivation of the scalar equation.  In this regard, in the next simulation, we will study the behavior of the  solution of the scalar equation starting from a different initial condition in which the numerical error is negligible compared to the modelling one. In particular, we will study its behavior in case $ \beta \gg0,$ considering $ \beta \approx0.8.$
    
    \subsection{Modeling error}
    Since we want measure the error due to the model, we have considered $\beta\gg 0.$ in this test. Therefore we set $A_g = 0.9$ and $\xi = 1/(1-\rho_0)$ where $\rho_0 = 0.2,$ and let us consider $u$ and $h$ such that $\beta \approx 0.84 $ and $F_r\approx 0.15.$ For this reason, the common settings of this test are: $[a,b]=[-2.5,20]$ the space domain; $N = 1600;$ $g = 9.81;$ $m_g=3;$ $t_{end} = 1000; $ $b(x)\equiv0$ and $$ u_0(x) = 1 + \delta e^{-\frac{(x+1)^2}{0.4^2}} $$ where $\delta = 0.006.$ $h$ and $z_b$ are keyed such that satisfied the scalar approximation with $h_0(a) = 4.21$ and $z_{b_0}(a) = 0.1.$ 
    
    Free boundary conditions are imposed at ghost points and CFL$_{scal} = 0.9$ while CFL$_{IMEX} = 1.$
   In Table \ref{Tab_mod_err}, we compute the relative error 
   $$ \rm{err} = \frac{||u_{scal} - u_{IMEX}||_1}{||u_{IMEX}||_1} $$
   between the solution computed with the scalar equation and semi-implicit numerical scheme. 
    \begin{table}[htbp]
        \begin{center}
            \begin{tabular}{|c|c|c|c|c|}
                \hline                        
                1600 points  & $h$ & $q$ & $z_b $ & $\eta$    \\ 
                \hline
                err      & 2.35E-3 &  1.13E-3 & 2.39E-3 & 2.14E-4 \\   
                \hline  
            \end{tabular}
            \end{center}
             \caption{ Errors in $L^1-$norm.}
             \label{Tab_mod_err}
    \end{table}
    As it can be observed, the relative error is of order $10^{-3}$ 
instead of the numerical discretization error that is of  order $10^{-6}$ and $10^{-7}$ respectively for both numerical schemes. For this reason, even if the scalar modeling gives a faster (in the computational sense) solution compared with the semi-implicit full system, the modeling error is not negligible when $\beta \gg0$ making it necessary the semi-implicit approach.

    Verified that the semi-implicit strategy leads to results similar to explicit and scalar approximations methods when $\beta \rightarrow0$ and confirmed that these results, in addition to being similar, are accurate with respect to the expected order, we want to explore the behavior and the results obtained in case a continuous waves group is imposed in the left boundary domain of the velocity $u.$

    \subsection{1D wave group}
    Let us consider the one-dimensional Exner system \eqref{Ex_sis_eta} and the second order semi-implicit method developed before. We want to verify, on the one hand, the temporal evolution of the sediment for very long times, for instance, a final time such that the initial dune has moved 10 times the initial amplitude; on the other hand, whether the presence of fast under-resolved surface waves has a significant effect in the evolution of the initial dune. In this way, we have three different time scales. The slowest one related to the velocity of the dune evolution $3\xi A_g u^2$ due to the Grass equation \eqref{q_b} if $m_g=3$ \cite{Grass}; the second one related to the water velocity $u;$ the fastest one related to the waves group of order $u+\sqrt{gh}.$
    \begin{figure}[!ht]
    	\centering	
    	\hspace{-1.2cm}
    	\includegraphics[height=0.15\textheight,width=1\textwidth]{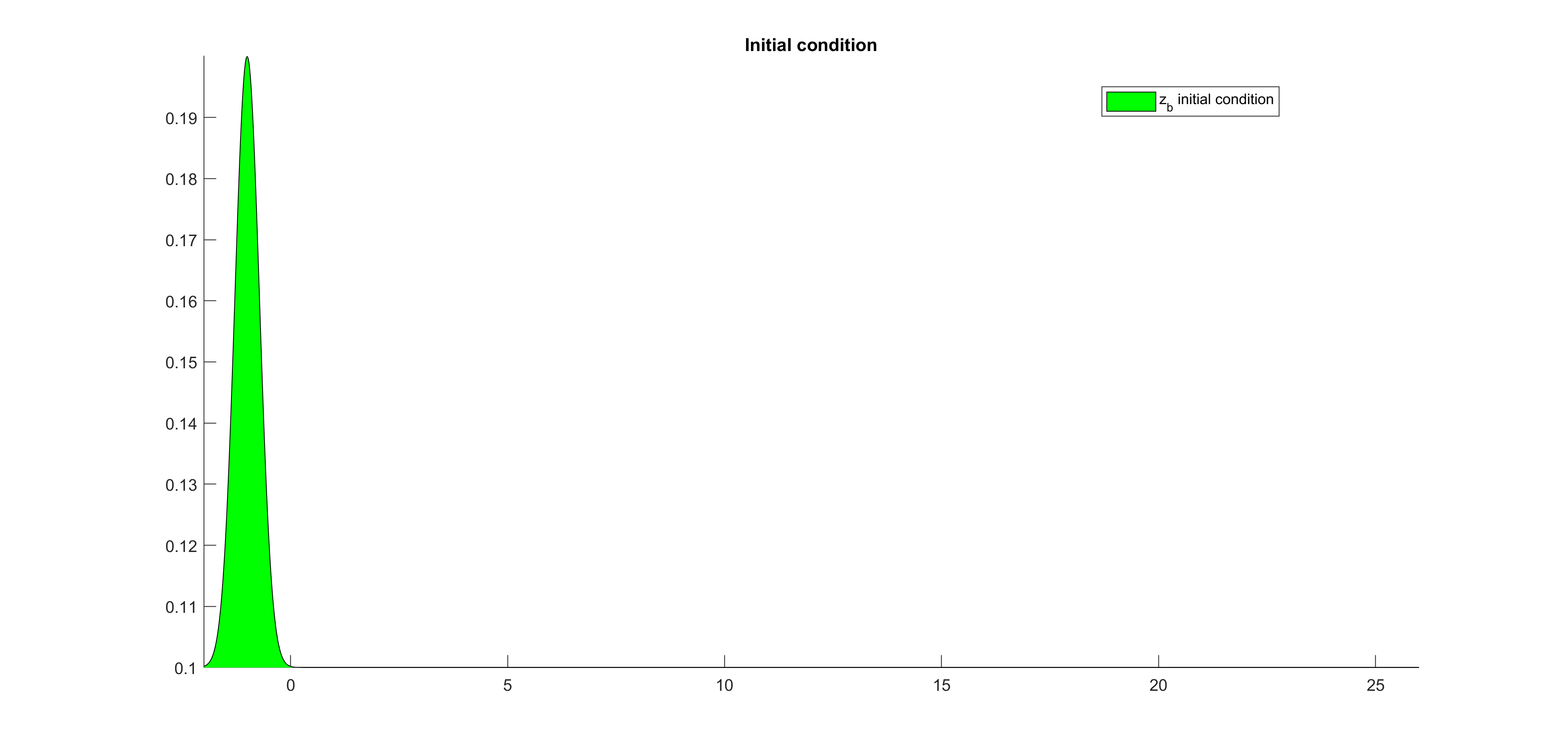}
    	\caption{Test \ref{Test_1}: (1D waves group). Initial condition of sediment for the Exner model on the interval $[-2,26]$ using a $2000-$mesh points.  }
    	\label{Test_1_1_waves}
    \end{figure}
    \begin{figure}[!ht]
    	\centering	
    	\hspace{-1.2cm}
    	\includegraphics[scale = 0.48]{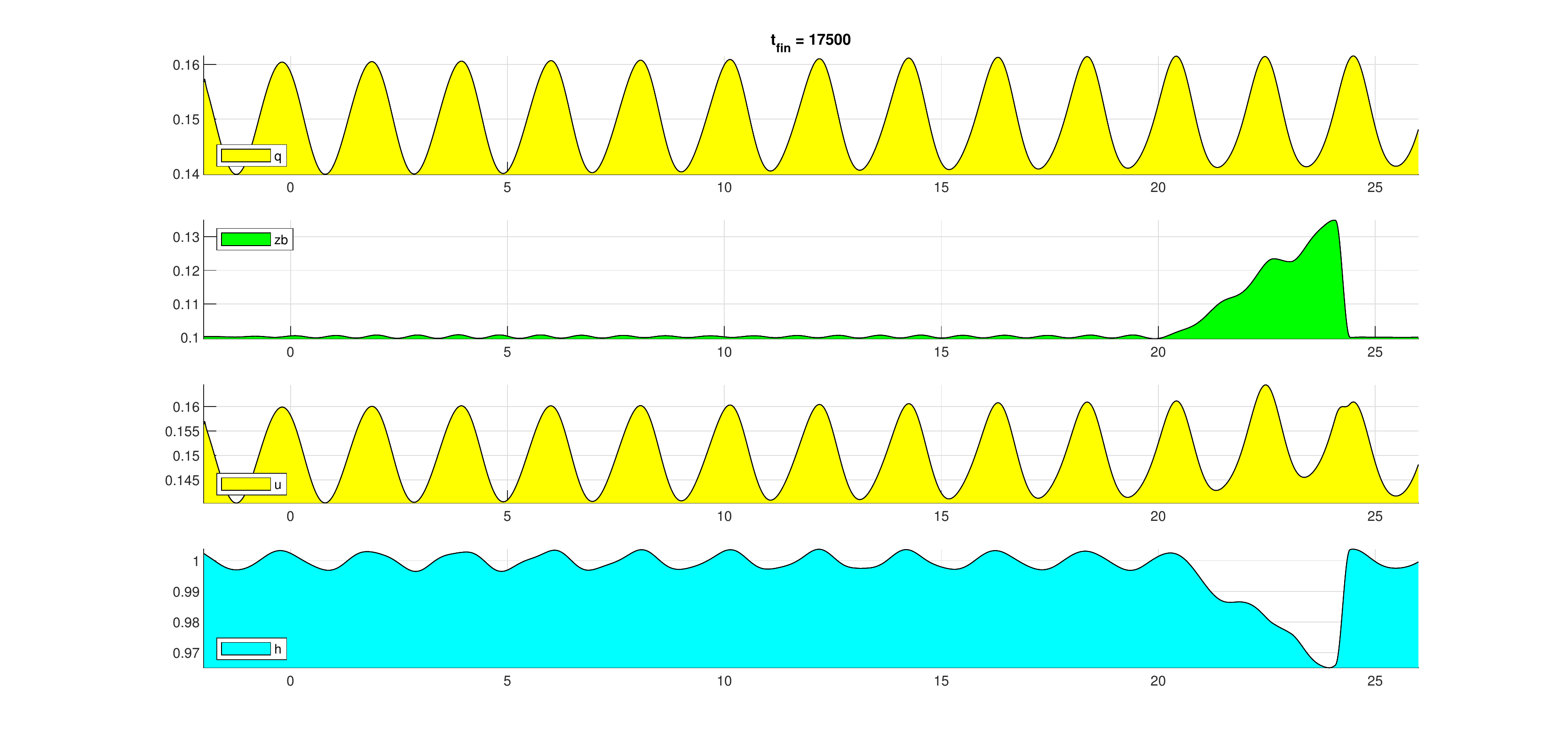}
    	\vspace{-0.8 cm}
    	\caption{Test \ref{Test_1}: (1D waves group). Numerical solutions of discharge (up), velocity (center-up), sediment layer (center-down) and thickness (down) for the Exner model on the interval $[-2,26]$ using a $2000-$mesh points at time $t=17500$ with CFL$ = 9.$}
    	\label{Test_1_2_waves}
    \end{figure}
    \begin{figure}[!ht]
    	\centering	
    	\hspace{-1.2cm}
    	\includegraphics[scale = 0.48]{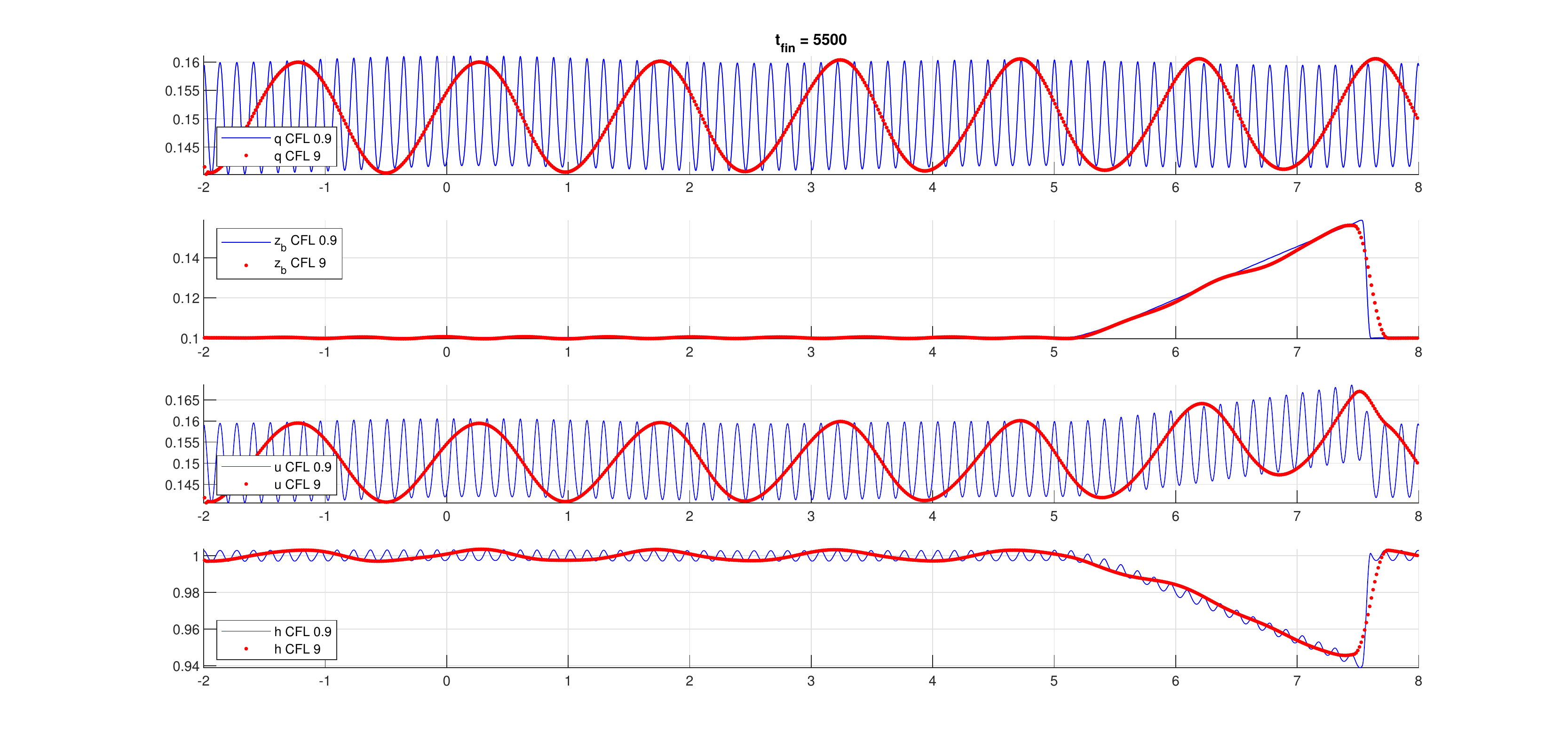}
    	\caption{Numerical solutions of discharge (up), velocity (center-up), sediment layer (center-down) and thickness (down) for the Exner model on interval $[-2,8]$ adopting a $2000-$mesh points at time $t = 5500$ with low CFL$=0.9$  and large CFL$=9.$}
    	\label{Test_1_3_waves}
    \end{figure}
    
    The settings of this test are: $[a,b]=[-2,26]$ the space domain; $A_g = 0.1;$ $\xi = 1/(1-\rho_0)$ where $\rho_0 = 0.2;$ $g = 9.81;$ $b(x) \equiv 0;$ $h_0(x) = 1;$ $u_0(x) = 0.15;$ \begin{equation}
        z_{b_0}(x)  = 0.1 + 0.1e^{-\frac{(x+1)^2}{0.4^2}}. 
    \end{equation}
    \begin{figure}[!ht]
    \vspace{-3cm}
    \centering
	\includegraphics[width=0.7 \textwidth]{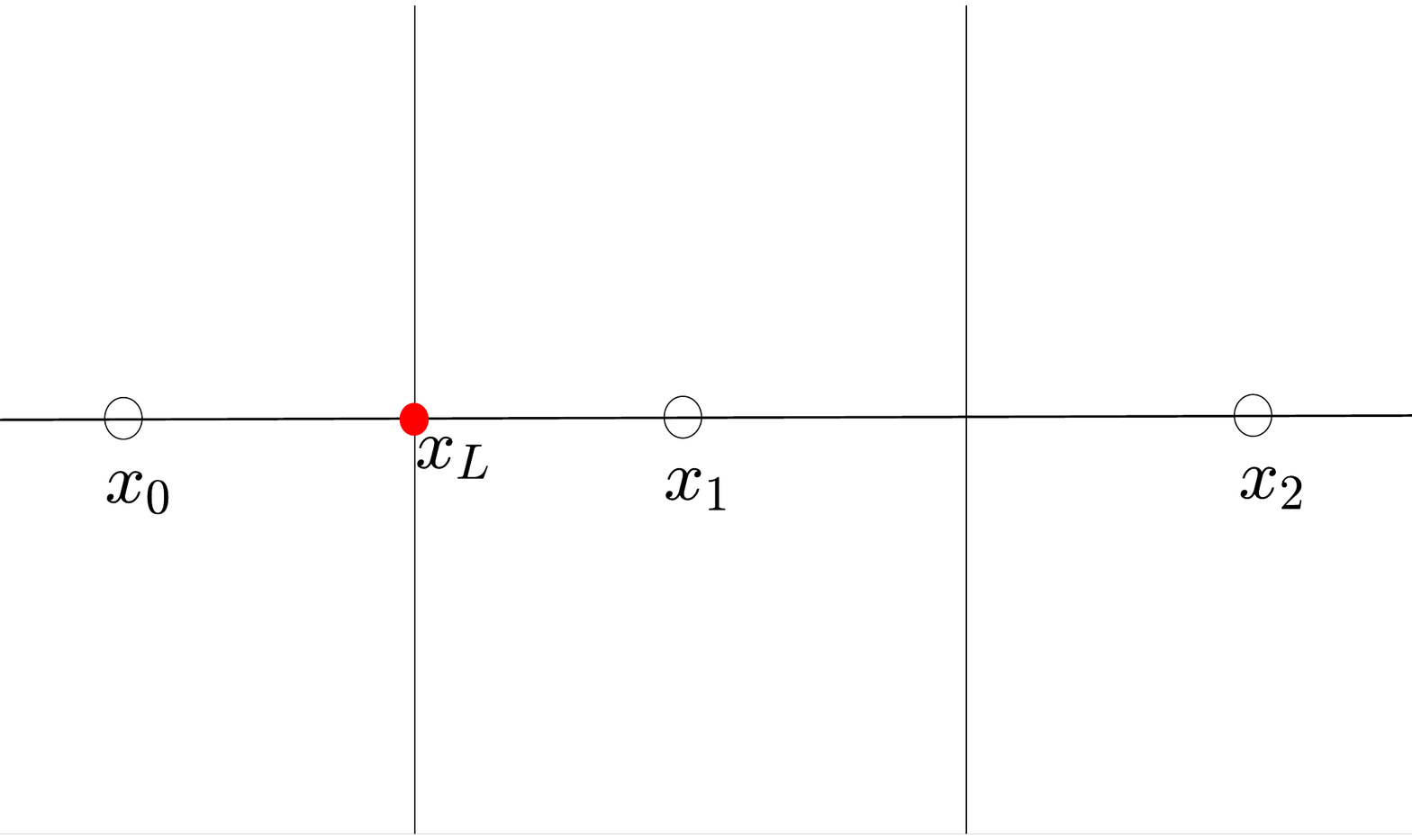}
	\vspace{-3cm}
	\caption{Left-boundary mesh with ghost point 0.}
	\label{Fig_bound}
\end{figure}
On the left boundary we impose the following conditions:
  $$
\begin{bmatrix}
      h_0 \\ u_0 
    \end{bmatrix} = \begin{bmatrix}
       h_1 + 1/g(\phi_t\Delta x + 0.5((u_1)^2 - \phi^2)) \\ 2\phi-u_1 
\end{bmatrix} 
$$
where $\phi$ and $\phi_t$ are respectively the waves group $\phi(t) = 0.15 + \mathcal{A}(sin(\omega t))$ and $\phi_t = d\phi/dt,$ in which $\mathcal{A}$ and $\omega$ are amplitude and frequency of the waves in our case set to 0.01 and 150 respectively \cite{TesiPhD,MaccaRussoBumi}.

Figures \ref{Test_1_1_waves}-\ref{Test_1_2_waves} show the initial and the numerical solutions for discharge, velocity, sediment layer and thickness obtained with the second-order semi-implicit scheme developed in the previous sections in which the stability condition is set  CFL$=9$ on the interval $[-2,26]$ adopting a $2000-$mesh points at time $t = 17500.$ 
    
Since we are imposing a very high-frequency signal to the left part of the domain, the water waves observed are not the real one. In fact, since the period of oscillations is $T = \frac{2\pi}{\omega} = 0.042$ and $\Delta t~0.038$ (with this settings), the ratio $\frac{T}{\Delta t} = 1.09$ which suggests that, more or less, at each time step a wave is inserted from the signal, so the visible waves on the graph are not the real waves but an understatement of them. To see clearly them a CFL reduction is necessary in order to have more time steps for each wave (see Figure \ref{Test_1_3_waves}). Nevertheless, we observe that, even if a shock appeared, the semi-implicit method with a low restriction in the stability condition (CFL$=9$), albeit diffusive, is able to capture and properly evolve the sedimentation. Furthermore, we can see how the surface waves group has a secondary effect on the sedimentation but still not negligible.

\section{2D Exner Model}
Let us consider the two-dimensional hyperbolic
shallow water equations
\begin{equation}
    \begin{cases}
        h_t + (hu)_x + (hv)_y = 0\\
        (hu)_t + (hu^2 + \ha gh^2)_x + (hvu)_y = -gh\frac{\partial b}{\partial x}\\
        (hv)_t + (huv)_x + (hv^2 + \ha gh^2)_y = -gh\frac{\partial b}{\partial y},
    \end{cases}
\end{equation}
where $(x,y)$ refers to the Cartesian plane $Oxy$ and $t$ is the time; $h(x,y,t),$ the thickness; $u(x,y,t)$ and $v(x,y,t),$ the horizontal velocity components; $g,$ the acceleration due to gravity; $b(x,y),$ the bottom topography. In particular, defining the momentum, $m = hu$ and $n = hv,$ we get:
    \begin{equation}
        \label{SH_2D}
        \begin{cases}
            h_t + (m)_x + (n)_y = 0\\
            (m)_t + (mu + \ha gh^2)_x + (mv)_y = -ghb_x\\
            (n)_t + (nu)_x + (nv + \ha gh^2)_y = -ghb_y.
        \end{cases}
    \end{equation}
The system of equations used in this section is obtained by coupling 2D shallow water equation \eqref{SH_2D} and the 2D sediment equation:
    \begin{equation}
        \label{sediment_equation_2D}
        (z_b)_t + (q_{x,b})_x + (q_{y,b})_y = 0
    \end{equation}
    where $z_b(x,y,t)$ represents the height of the sediment layer and, $q_{x,b}(u,v)$ and $q_{y,b}(u,v),$ the solid transport discharge parameters, in our case computed by the 2D Grass model \cite{Grass,Hudson,MURILLO20108704} 
    \begin{align}
        \label{q_b_x}
        q_{x,b} =& \xi A_gu(u^2+v^2)^{\frac{m_g-1}{2}} \\
        \label{q_b_y}
        q_{y,b} =& \xi A_gv(u^2+v^2)^{\frac{m_g-1}{2}},
    \end{align}
    with $m_g\in[1,4],$ $A_g\in ]0,1[$ and $\xi = 1/(1-\rho_0)$ where $\rho_0$ is the porosity of the sediment layer.
    
    In this way, the 2D Exner system of balance laws is given by:
    \begin{equation}
        \label{2D_Ex_sis}
        \begin{cases}
             h_t + (m)_x + (n)_y = 0\\
            (m)_t + (mu + \ha gh^2)_x + (mv)_y = -gh(b + z_b)_x\\
            (n)_t + (nu)_x + (nv + \ha gh^2)_y = -gh(b + z_b)_y\\
            (z_b)_t + (q_{x,b})_x + (q_{y,b})_y = 0.
        \end{cases}
    \end{equation}
Following the same procedure done for the 1D case, system \eqref{2D_Ex_sis} can be written in the following way:
\begin{equation}
    \label{2D_Hyp_sys} 
    \partial_tU + A_1(U)\partial_xU + A_2\partial_yU =0 
    \end{equation}
where
\[
    U = \begin{bmatrix}
            h \\ m \\ n \\ z_b
       \end{bmatrix}; 
    \quad A_1(U) = \begin{bmatrix}
            0 & 1 & 0 & 0 \\
            gh-u^2 & 2u & 0 & gh \\
            -uv & v & u & 0 \\
            \alpha_x & \beta_x & \gamma_x & 0
        \end{bmatrix}; \quad A_2(U) = \begin{bmatrix}
            0 & 0 & 1 & 0 \\
            -uv & v & u & 0 \\
            gh-v^2 & 0 & 2v & gh \\
            \alpha_y & \beta_y & \gamma_y & 0
        \end{bmatrix};
\]
in which the terms $\alpha_s, \beta_s$ and $\gamma_s,$ with $s\in \{x,y\},$ represent the $q_{s,b}$ derivatives respect to $h,m$ and $n,$ i.e. for $s\in \{x,y\}$
\[
    \alpha_s = \frac{\partial q_{s,b}}{\partial h},\quad \beta_s = \frac{\partial q_{s,b}}{\partial m},\quad  \gamma_s = \frac{\partial q_{s,b}}{\partial n}.
\]
In particular, defining $c_1 = \Bigl(\xi A_g (u^2 + v^2)^{\frac{m_g-1}{2}}\Bigr)/h $ and $c_2 = \Bigl(\xi A_g (u^2 + v^2)^{\frac{m_g-3}{2}}\Bigr)/h $ we obtain
\begin{align*}
    \alpha_x &= -(m_g-1)c_1u \\ 
    \alpha_y &= -(m_g-1)c_1v \\ 
    \beta_x &= c_1 + (m_g-1)c_2u^2 \\
    \beta_y & = (m_g-1)c_2uv \\
    \gamma_x &= \beta_y \quad {\rm{and}} \quad \gamma_y = \beta_x. 
\end{align*}

At the end, let us rewrite the 2D Exner system \eqref{2D_Ex_sis} in function of $\eta$ where $\eta(x,y,t) = h(x,y,t) + b(x,y) + z_b(x,y,t)$ represents the elevation of the undisturbed water surface. In practise, system \eqref{2D_Ex_sis} becomes:
    \begin{equation}
        \label{2D_Ex_sis_eta}
        \begin{cases}
             \eta_t + (m + q_{x,b})_x + (n + q_{y,b} )_y = 0\\
            (m)_t + (mu)_x + (mv)_y + gh(\eta)_x = 0\\
            (n)_t + (nu)_x + (nv)_y + gh(\eta)_y = 0\\
            (z_b)_t + (q_{x,b})_x + (q_{y,b})_y = 0.
        \end{cases}
    \end{equation}

\section{2D semi implicit scheme} 
\label{2D_semi_impl_sect}
In this section, we present the extension of the  semi implicit scheme to the 2D Exner model \eqref{2D_Ex_sis_eta}. We plan to treat implicitly the surface water waves while the corresponding slow sediment wave will be treated explicitly. As we have done in Section \ref{1D semi implicit}, we will just consider first and second order schemes in space and time.
    
We consider a rectangular computational domain $\Omega \equiv [a_1,b_1]\times[a_2,b_2]$, and divide it into $N_x\times N_y$ cells defined by $ I_{i,j} = [x_{i-\ha},x_{i+\ha}]\times [y_{j-\ha},y_{j+\ha}],$ with $i = 1,\ldots,N_x$ and $j = 1,\ldots,N_y$.  For the sake of simplicity, we adopt a uniform Cartesian mesh direction by direction with mesh spacing, respectively, $\Delta x$ and $\Delta y,$ i.e. $x_i = a_1 + (i-\ha)\Delta x$ and $y_j = a_2 + (j-\ha)\Delta y.$ As previously, $\Delta t$ is the time step such that $t^k = k\Delta t.$\footnote{In practice it is better to assign $\Delta t$  dynamically at each time step by imposing some CFL condition. The choice of constant time step here is adopted in order to simplify the notation in the description of the method.}
    
Finally, we denote by $U_{i,j}^k$ an approximation on the mean value of $U$ over cell $I_{i,j}$ at time $t = t^k$ as: $$U_{i,j}^k \cong \frac{1}{\Delta x\Delta y}\int_{x_{i-\ha}}^{x_{i+\ha}} \int_{y_{j-\ha}}^{y_{j+\ha}} U(x,y,t^k)\,dy\,dx.$$

\subsection{First order scheme}
Let us consider the system \eqref{2D_Ex_sis_eta}, following the idea probosed by \cite{Boscarino-Filbet} in which IMEX Runge-Kutta methods are used for systems in which the stiffness is not necessarily of additive or partitioned type. Let us write system \eqref{2D_Ex_sis_eta} as a large system of ODE’s, in which we adopt suitable discrete operators for the approximation of space derivatives. The key point in \cite{Boscarino-Filbet} is to identify which specific term has to be treated implicitly and which can be treated explicitly. With this in mind, let us write \eqref{2D_Ex_sis_eta} in form 
\begin{equation}
        \label{2D_ode_bosca_semi}
        U' = H(U_E,U_I).
    \end{equation}
    where $U = [\eta, m,n, z_b]^T$ and $H(U_E,U_I)$ is given by
    \begin{equation}
        \label{2D_ode_form_semi}
        H(U_E,U_I) = 
        \begin{bmatrix}
            -(m_I +(q_{x,b})_E)_x - (n_I +(q_{y,b})_E)_y \\
            -( (mu)_E)_x + (mv)_E)_y - gh_E(\eta_I)_x \\
            -( (nu)_E)_x + (nv)_E)_y - gh_E(\eta_I)_y \\
            -(q_{x,b})_E)_x - (q_{y,b})_E)_y
        \end{bmatrix}
    \end{equation}
where the subscript $E$ and $I$ denote which term has to be treated explicitly and which implicitly.  

The semi-implicit scheme can be written in the form
    \begin{equation}
        \label{2D_ode_bosca}
        U' = \tilde{H}(U_E,U_I).
    \end{equation}
with
    \begin{equation}
        \label{2D_ode_form_part}
        \tilde{H}(U_E,U_I) = 
        \begin{bmatrix}
           -\hat{D}_x((q_{x,b})_E) - \hat{D}_y((q_{y,b})_E)  & -D_x(m_I) -D_y(n_I) \\ -\hat{D}_x((mu)_E) -\hat{D}_y((mv)_E) &  -gh_ED_x(\eta_I) \\
           -\hat{D}_x((nu)_E) -\hat{D}_y((nv)_E) &  -gh_ED_y(\eta_I) \\-\hat{D}_x((q_{x,b})_E) - \hat{D}_y((q_{y,b})_E) &
        \end{bmatrix}.
    \end{equation}

A semi-discrete in time first order semi implicit method can be written as:
\begin{equation}
    \begin{cases}
        \label{2D First order}
        \eta^{k+1} = \eta^k - \Delta t \hat{D}_x(q_{x,b}^k) - \Delta t \hat{D}_y(q_{y,b}^k) - \Delta t D_x(m^{k+1}) - \Delta t D_y(n^{k+1}),\\
        m^{k+1} = m^k - \Delta t \hat{D}_x(m^ku^k) - \Delta t\hat{D}_y(m^kv^k) - \Delta tgh^kD_x(\eta^{k+1}), \\
         n^{k+1} = n^k - \Delta t \hat{D}_x(n^ku^k) - \Delta t\hat{D}_y(n^kv^k) - \Delta tgh^kD_y(\eta^{k+1}),\\
        z_b^{k+1} = z_b^k - \Delta t \hat{D}_x(q_{x,b}^k) - \Delta t \hat{D}_y(q_{y,b}^k),
    \end{cases}
\end{equation}
where the differential operators $D_x,D_y,\hat{D}_x$ and $\hat{D}_y$ are defined as in Section \ref{1D semi implicit}, direction by direction. In particular:
\begin{itemize}
	    \item $D_x(\tilde{F}_{i,j}) = \frac{\tilde{F}_{i+\ha,j} - \tilde{F}_{i-\ha,j}}{\Delta x},$ in which $\tilde{F}_{i\pm\ha,j}$ are suitably defined on cell edges; 
	    \item $D_y(\tilde{G}_{i,j}) = \frac{\tilde{G}_{i,j+\ha} - \tilde{G}_{i,j-\ha}}{\Delta y},$ in which $\tilde{G}_{i,j\pm\ha}$ are suitably defined on cell edges; 
	    \item $\hat{D}_x(\tilde{F}_{i,j}) = \frac{\tilde{F}_{i+\ha,j} - \tilde{F}_{i-\ha,j}}{\Delta x},$ where  $\tilde{F}_{i+\ha,j} = \ha\Bigl( \tilde{F}(U_{i+\ha,j}^{-}) + \tilde{F}(U_{i+\ha,j}^{+}) - \alpha_{x,i+\ha,j}\big(U_{i+\ha,j}^{+} - U_{i+\ha,j}^{-}\big)\Bigr)$ is the Rusanov flux and $\alpha_{x,i+\ha,j}$ is related to the eigenvalues of the explicit sub system. In our case,  $\alpha_{x} \approx |u|;$
	    \item $\hat{D}_y(\tilde{G}_{i,j}) = \frac{\tilde{G}_{i,j+\ha} - \tilde{G}_{i,j-\ha}}{\Delta y},$ where  $\tilde{G}_{i,j+\ha} = \ha\Bigl( \tilde{G}(U_{i,j+\ha}^{-}) + \tilde{G}(U_{i,j+\ha}^{+}) - \alpha_{y,i,j+\ha}\big(U_{i,j+\ha}^{+} - U_{i,j+\ha}^{-}\big)\Bigr)$ is the Rusanov flux and $\alpha_{y,i,j+\ha}$ is related to the eigenvalues of the explicit sub system. In our case,  $\alpha_{y} \approx |v|.$
	\end{itemize}

To emphasize the explicit part from the implicit one, let us rewrite system \eqref{2D First order} as:
\begin{equation}
    \label{2D First order star}
	\begin{cases}
        m^* = m^k - \Delta t\hat{D}_x(m^ku^k) -  \Delta t\hat{D}_y(m^kv^k); \\
        n^* = n^k - \Delta t\hat{D}_x(n^ku^k) -  \Delta t\hat{D}_y(n^kv^k); \\
        \eta^* = \eta - \Delta t\hat{D}_x(q_{x,b}^k) - \Delta t\hat{D}_y(q_{y,b}^k) - \Delta t\hat{D}_x(m^*) - \Delta t\hat{D}_y(n^*) ;\\
        \eta^{k+1} = \eta^* + g \Delta t^2D_x(h^kD_x(\eta^{k+1})) + g \Delta t^2D_y(h^kD_y(\eta^{k+1})); \\
        m^{k+1} = m^* - g\Delta th^kD_x(\eta^{k+1});\\
        n^{k+1} = n^* - g\Delta th^kD_y(\eta^{k+1});\\ z_b^{k+1} = z_b^k - \Delta t\hat{D}_x(q_{x,b}^k) - \Delta t\hat{D}_y(q_{y,b}^k).
	\end{cases}
\end{equation}
The procedure to solve system \eqref{2D First order}, hence \eqref{2D First order star}, is:
\begin{enumerate}
    \item compute $m^* = m^k - \Delta     t\hat{D}_x(m^ku^k) -  \Delta t\hat{D}_y(m^kv^k)$ as 
    $$ m^*_{i,j} = m^k_{i,j} - \frac{\Delta t}{\Delta x}\Bigl(F^k_{i+\ha,j,m} - F^k_{i-\ha,j,m}\Bigr) - \frac{\Delta t}{\Delta y}\Bigl(G^k_{i,j+\ha,m} - G^k_{i,j-\ha,m}\Bigr), $$ where $F_m$ and $G_m$ are computed direction by direction with the corresponding Rusanov flux defined \modify{in} Section \ref{2D_semi_impl_sect} related to $mu$ and $mv;$ 
    \item compute $n^* = n^k - \Delta t\hat{D}_x(n^ku^k) -  \Delta t\hat{D}_y(n^kv^k)$ as 
    $$ n^*_{i,j} = n^k_{i,j} - \frac{\Delta t}{\Delta x}\Bigl(F^k_{i+\ha,j,n} - F^k_{i-\ha,j,n}\Bigr) - \frac{\Delta t}{\Delta y}\Bigl(G^k_{i,j+\ha,n} - G^k_{i,j-\ha,n}\Bigr), $$
    where $F_n$ and $G_n$ are computed direction by direction with the Rusanov flux related to $nu$ and $nv;$\footnote{In this case there is an abuse of notation, the $n$ at the subscript indicates the momentum in the $y$ direction while the $n$ at the apex the generic temporal level.}
     \item compute $\eta^* = \eta - \Delta t\hat{D}_x(q_{x,b}^k) - \Delta t\hat{D}_y(q_{y,b}^k) - \Delta t\hat{D}_x(m^*) - \Delta t\hat{D}_y(n^*) $ as
    \begin{align*}
        \eta^*_{i,j} = \eta_{i,j} - &\frac{\Delta t}{\Delta x} \Bigl(F^k_{i+\ha,j\eta} - F^k_{i-\ha,j,\eta}\Bigr) -  \frac{\Delta t}{\Delta y} \Bigl(G^k_{i+\ha,j,\eta} - G^k_{i-\ha,j,\eta}\Bigr)+\\ -&\frac{\Delta t}{2\Delta x} \Bigl(m^*_{i+1,j} - m^*_{i-1,j}\Bigr) -  \frac{\Delta t}{2\Delta y} \Bigl(n^*_{i,j+1} - n^*_{i,j-1}\Bigr),
    \end{align*}
    where $F^3$ and $G^3$ are the Rusanov operators referred to $q_{x,b}$ and $q_{y,b}.$
    \item Let be $k_x \equiv g({\Delta t}/{\Delta x})^2$ and $k_y \equiv g({\Delta t}/{\Delta y})^2.$ Solve the equation for $\eta^{n+1}$ 
    \[
        \eta^{k+1} = \eta^* + g \Delta t^2D_x(h^kD_x(\eta^{k+1})) + g \Delta t^2D_y(h^kD_y(\eta^{k+1}))
    \]
    which can be written in components as
    \begin{align*}
        &\eta_{i,j}^{k+1}\Bigl(1-k_x(h_{i,j-\ha}^k + h_{i,j+\ha}^k) - k_y(h_{i-\ha,j}^k + h_{i+\ha,j}^k)\Bigr) +\\
        +&\eta_{i,j-1}^{k+1}\Bigl(k_x(h_{i,j-\ha}^k\Bigr) + \eta_{i,j+1}^{k+1}\Bigl(k_y(h_{i,j+\ha}^k\Bigr) +\\
        +&\eta_{i-1,j}^{k+1}\Bigl(k_x(h_{i-\ha,j}^k\Bigr) + \eta_{i+1,j}^{k+1}\Bigl(k_y(h_{i+\ha,j}^k\Bigr) = \eta^*_{i,j},
    \end{align*}
    where $h_{i\pm \ha,j}^k = \ha(h_{i,j}^k + h_{i\pm 1,j}^k)$ and $h_{i,j\pm \ha}^k = \ha(h_{i,j}^k + h_{i,j\pm +1}^k).$
    
    This is an invertible linear system which can be solved to detect $\eta^{k+1} = [\eta^{k+1}_{i,j}]$ for all $i = 1,\ldots,N_x$ and $j=1,\ldots,N_y;$
    \item compute $m^{k+1} = m^* - g\Delta th^kD_x(\eta^{k+1})$ as 
    $$ m^{k+1}_{i,j} = m^*_{i,j} - g\frac{\Delta t}{\Delta x}h_{i,j}^k\Bigl(\eta_{i+\ha,j}^{k+1} - \eta_{i-\ha,j}^{k+1}\Bigr),$$ where $\eta_{i\pm\ha,j}^{k+1} = \ha\Bigl(\eta_{i,j}^{k+1}+\eta_{i\pm 1,j}^{k+1}\Bigr)$ for all $j = 1,\ldots,N_y;$ $i=1...N_x.$ 
    \item compute $n^{k+1} = n^* - g\Delta th^kD_y(\eta^{k+1})$ as 
    $$ n^{k+1}_{i,j} = n^*_{i,j} - g\frac{\Delta t}{\Delta y}h_{i,j}^k\Bigl(\eta_{i,j+\ha}^{k+1} - \eta_{i,j-\ha}^{k+1}\Bigr),$$ where $\eta_{i,j\pm\ha}^{k+1} = \ha\Bigl(\eta_{i,j}^{k+1}+\eta_{i,j\pm1}^{k+1}\Bigr)$ for all $i = 1,\ldots,N_x;$ $j=1...N_y;$
    \item compute $z_b^{k+1} = z_b^k - \Delta t\hat{D}_x(q_{x,b}^k) - \Delta t\hat{D}_y(q_{y,b}^k)$ as
       $$ z_{b_{i,j}}^{k+1} = z_{b_{i,j}}^{k+1} - \frac{\Delta t}{\Delta x} \Bigl(F^k_{i+\ha,j,z_b} - F^k_{i-\ha,j,z_b}\Bigr) -  \frac{\Delta t}{\Delta y} \Bigl(G^k_{i+\ha,j,z_b} - G^k_{i-\ha,j,z_b}\Bigr),$$
       in which $F_{z_b}$ and $G_{z_b}$ are computed direction by direction with the Rusanov flux related to $q_{x,b}$ and $q_{y,b};$
   \item compute $ h_{i,j}^{k+1} = \eta_{i,j}^{k+1} - b_{i,j} - z_{b_{i,j}}^{k+1},$ for all $i = 1,\ldots,N_x$ and for all $j = 1,\ldots,N_y.$
\end{enumerate}

	\subsection{Second order scheme}\label{2D_Sec_ord_scheme}
As have been done for the first order case in 1D and following \cite{Boscarino-Filbet}, we write system \eqref{2D_Ex_sis_eta} in the ODE form \eqref{ode_bosca_semi}-\eqref{ode_form_semi}   

After space discretization, the semi-implicit scheme can be written in the form
	\begin{equation}
	    \label{2D ode_bosca}
	    U' = H(U,U),
	\end{equation}
	where, in the 2D case, $U=[\eta,m,n,z_b]^T$ and $H(U,U)$ is defined as:
	\begin{equation}
	    \label{2D ode_form}
	    H(U,U) = \begin{bmatrix}
	       -( q_{x,b} + m)_x - ( q_{y,b} + n)_y \\ 
	       -( mu)_x - gh(\eta)_x - (mv)_y \\
	       -(nu)_x - (nv)_y -gh(\eta)_y \\
	       -(q_{x,b})_x -(q_{y,b})_y 
	    \end{bmatrix}
	\end{equation}
	that, differentiating between explicit and implicit part, the system is:
    \begin{equation}
	    \label{2D ode_form_partitioned}
	    H(U_E,U_I) = \begin{bmatrix}
	       -\hat{D}_x( (q_{x,b})_E) - \hat{D}_y( (q_{y,b})_E) &  - D_x(m_I) - D_y(n_I) \\ 
	       -\hat{D}_x((mu)_E) -\hat{D}_y((mv)_E) & - gh_ED_x(\eta_I) \\
	       -\hat{D}_x((nu)_E) - \hat{D}_y((nv)_E) &  -gh_ED_y(\eta_I) \\
	       -\hat{D}_x((q_{x,b})_E) - \hat{D}_y((q_{y,b})_E) & 
	    \end{bmatrix}.
	\end{equation}
    With this in mind, we apply an IMEX scheme to system \eqref{2D ode_form_partitioned}.
    The general procedure to update the numerical solution from time $t_k$ to $t_{k+1}$ using an $s$-stage Runge-Kutta IMEX method is the following:
    \begin{itemize}
        \item Stage values: For $i=1,\ldots,s$ compute
        \begin{align*}
            U^{(i)}_E & = U^k + \Delta t\sum_{j=1}^{i-1}a_{i,j}^EH\left(U_E^{(j)},U_I^{(j)}\right)\\
            U^{(i)}_I & = U^k + \Delta t
            \left(\sum_{j=1}^{i-1}a_{i,j}^I H\left(U_E^{(j)},U_I^{(j)}\right) + a_{i,i}^I H\left(U_E^{(i)},U_I^{(i)}\right)\right).
        \end{align*}
        \item Numerical solution:\\
        $U^{k+1} = U_I^{(s)}.$
    \end{itemize}    
    
Here we consider the IMEX scheme defined by the following  double Butcher tableau \cite{Boscarino-Filbet}: 
\begin{equation}
\label{tableau_2D}
    \begin{array}{c|cc}
             & 0 &  \\
            c & c & 0\\ \hline
            & 1-\gamma & \gamma
        \end{array} 
        \hspace{3 cm}
        \begin{array}{c|cc}
            \gamma & \gamma &  \\
            1 & 1-\gamma & \gamma\\ \hline
            & 1-\gamma & \gamma
        \end{array}
\end{equation}
where $\gamma = 1 - \frac{1}{\sqrt{2}}$ and $c = \frac{1}{2\gamma}.$

	Following the same reconstruction used for the 1D model, the procedure to update the numerical solution for \eqref{2D ode_bosca} is:
	\begin{enumerate}
	    \item $U_E^{(1)} = U^k;$
        \item $U_I^{(1)} = U^k + \Delta t\gamma H(U_E^{(1)},U_I^{(1)});$
        \item $U_E^{(2)} = (1-\frac{c}{\gamma})U^k + \frac{c}{\gamma}U_I^{(1)};$
        \item $U_I^{(2)} = (1-\frac{1-\gamma}{\gamma})U^k + \frac{1-\gamma}{\gamma}U_I^{(1)} +\Delta t\gamma H(U_E^{(2)},U_I^{(2)});$
        \item $U^{k+1} = U_I^{(2)}.$
	\end{enumerate}

    \section{2D Exner numerical experiments} \label{Test_2}
    In this section we test the semi-implicit scheme for the 2D Exner model with two different initial conditions: a parabolic and a conical sediment. As boundary conditions, zero Neumann conditions are imposed at ghost cells. As common settings of these experiments: $A_g = 0.1;$ $m_g = 3$ and $\rho_0 = 0.2.$

    \subsection{Parabolic Sediment}\label{Test_2_a}
\begin{figure}[!ht]
    	\centering	
    	\hspace{-1.2cm}
    	\includegraphics[scale = 0.5]{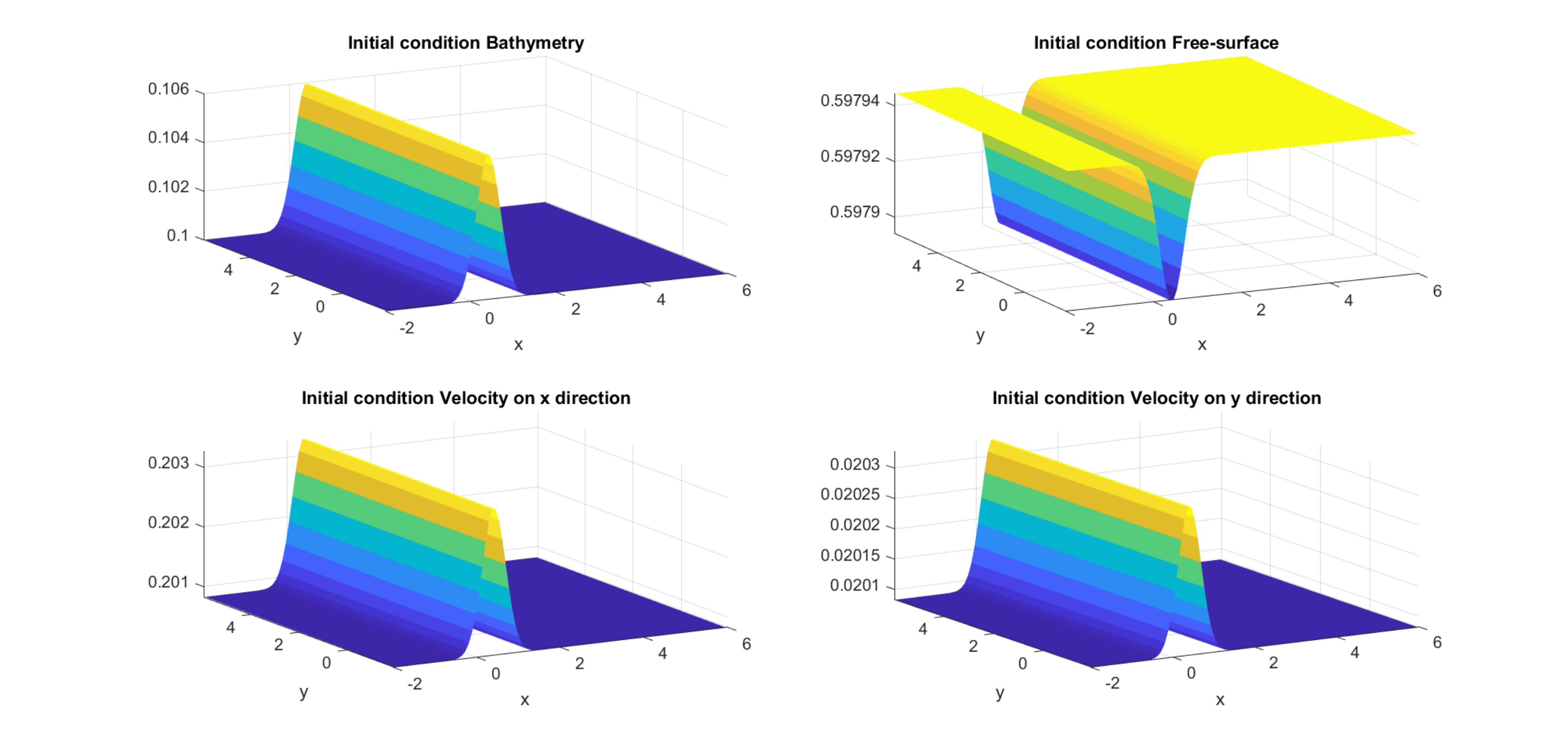}
    	\vspace{-0.8 cm}
    	\caption{Test \ref{Test_2_a}: (2D Exner parabolic sediment). Initial condition of sediment layer (top-left), thickness (top-right) and velocity (down) for the 2D Exner model on the square $[-2,6]\times[-2,6]$ using a $100\times100$ mesh points. }
    	\label{Test_2_1}
    \end{figure}

    With this purpose in mind, let us consider the 2D Exner model \eqref{2D_Ex_sis_eta} where initial conditions are so set: $\eta_0(x,y,0) = 0.6,$ $b(x,y) = 0,$ $m(x,y,0) = 0.1,$ $n(x,y,0)=0.01$\footnote{A small value is imposed in the flow in the $y-$direction to have a purely two-dimensional example. If this initial parameter is set to zero we will have a contribution only in the $x-$direction thus making the problem one dimensional.} 
   and 
    \begin{equation}
        \label{zb0} 
        z_{b_0}(x,y) = 0.1 + 0.006e^{-\frac{(x-0.4)^2}{0.4^2}}
    \end{equation}
    a one-dimensional parabolic sediment, see Figure \ref{Test_2_1}.    
    \begin{figure}[!ht]
    	\centering	
    	\hspace{-1.2cm}
    	\includegraphics[scale = 0.5]{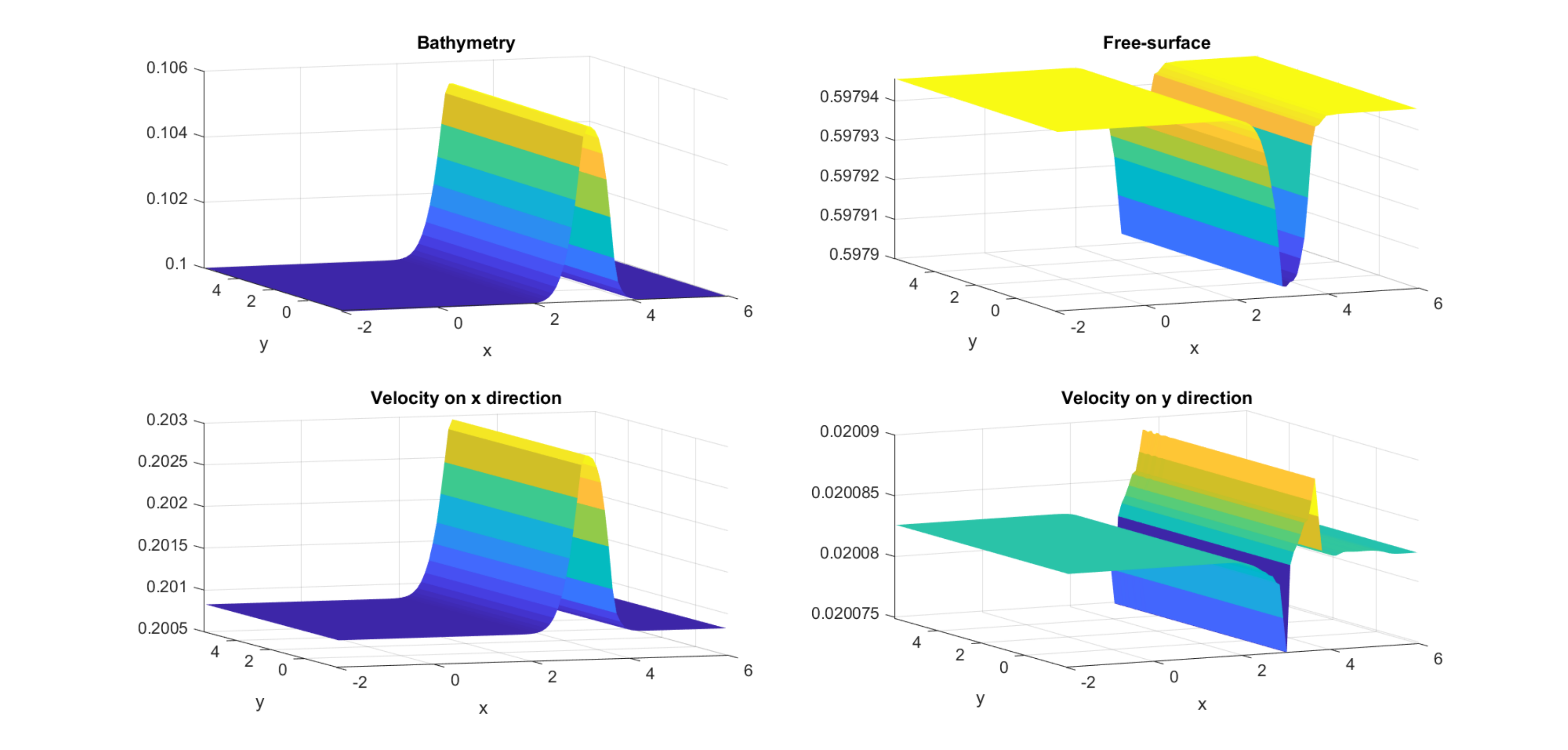}
    	\vspace{-0.8 cm}
    	\caption{Test \ref{Test_2_a}: (2D Exner parabolic sediment). Numerical solution for sediment layer (top-left), thickness (top-right) and velocity (down) for the 2D Exner model on the square $[-2,6]\times[-2,6]$ using a $100\times100$ mesh points at time $t=450$ with CFL$ = 6.$  }
    	\label{Test_2_2}
    \end{figure}
    The numerical results are obtained with the second order semi-implicit scheme introduced on Section \ref{2D_Sec_ord_scheme} on the square $[-2,6]\times[-2,6]$ adopting a $100\times100$ mesh points, CFL$=9$ at time $t = 450.$ As it shown in Figure \ref{Test_2_2}, the numerical results are in accordance with the one-dimensional one.

    \subsection{Conical Sediment} \label{Test_2_b}
\begin{figure}[!ht]
    	\centering	
    	\hspace{-1.2cm}
    	\includegraphics[scale = 0.5]{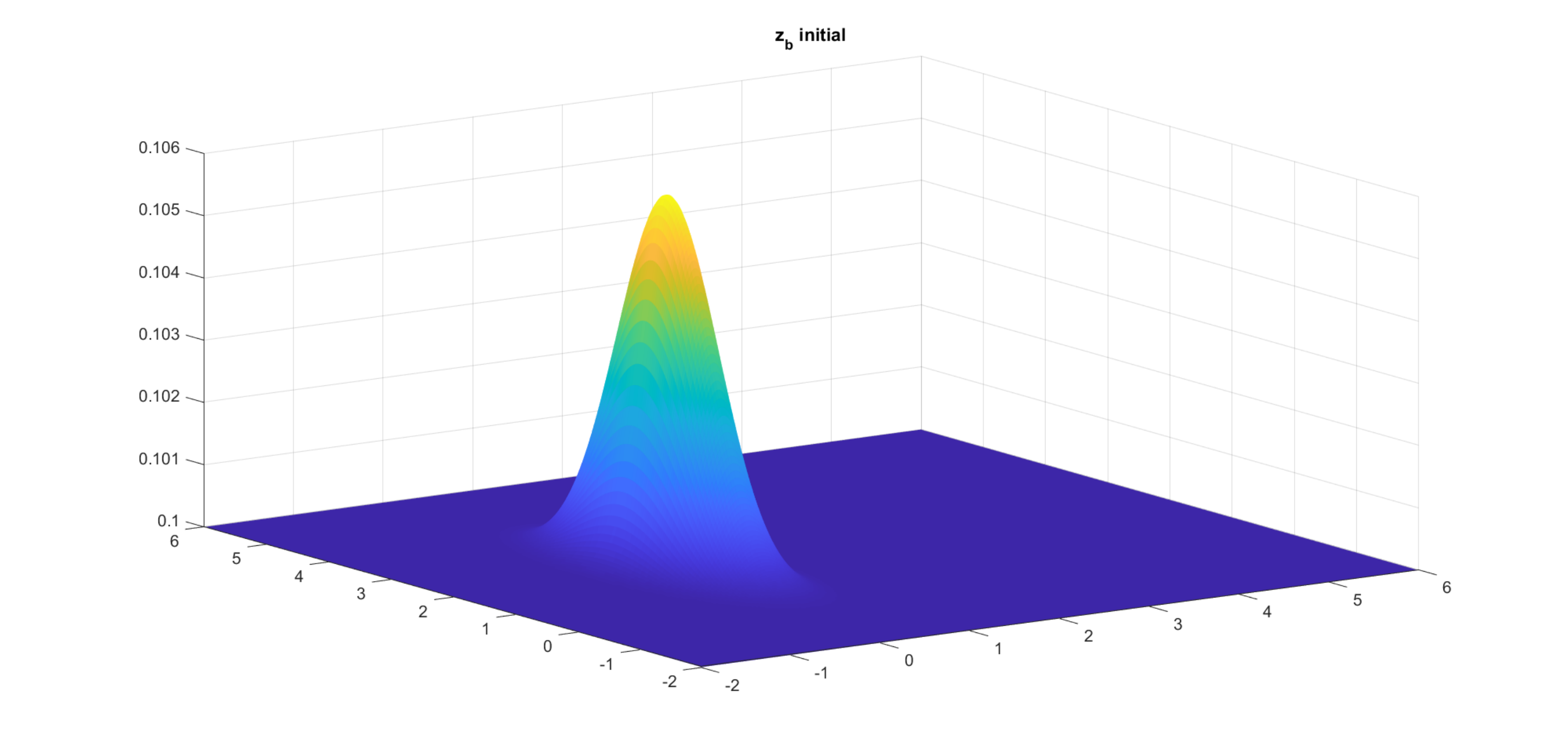}
    	\vspace{-0.8 cm}
    	\caption{Test \ref{Test_2_b}: (2D Exner conical sediment). Initial condition of sediment layer for the 2D Exner model on the square $[-2,6]\times[-2,6]$ using a $300\times300$ mesh points. }
    	\label{Test_2_1_con}
    \end{figure}

    As last class of experiments we test a two-dimensional conical sediment with a zero-flow in the $y-$direction ($n(x,y,0)=0$ and a complete 2D cases. As first one, let us consider the 2D Exner model \eqref{2D_Ex_sis_eta} where initial conditions are so set: $\eta_0(x,y,0) = 1.8,$ $b(x,y) = 0,$ $m(x,y,0) = 0.3,$ $n(x,y,0)=0$ and 
    \begin{equation}
        \label{zb0_1} 
        z_{b_0}(x,y) = 0.1 + 0.006e^{-\frac{(x-0.4)^2}{0.4^2} - (y-3)^2}
    \end{equation}
    a conical sediment, see Figure \ref{Test_2_1_con}.
    
    \begin{figure}[!ht]
    	\centering	
    	\hspace{-1.2cm}
    	\includegraphics[scale = 0.5]{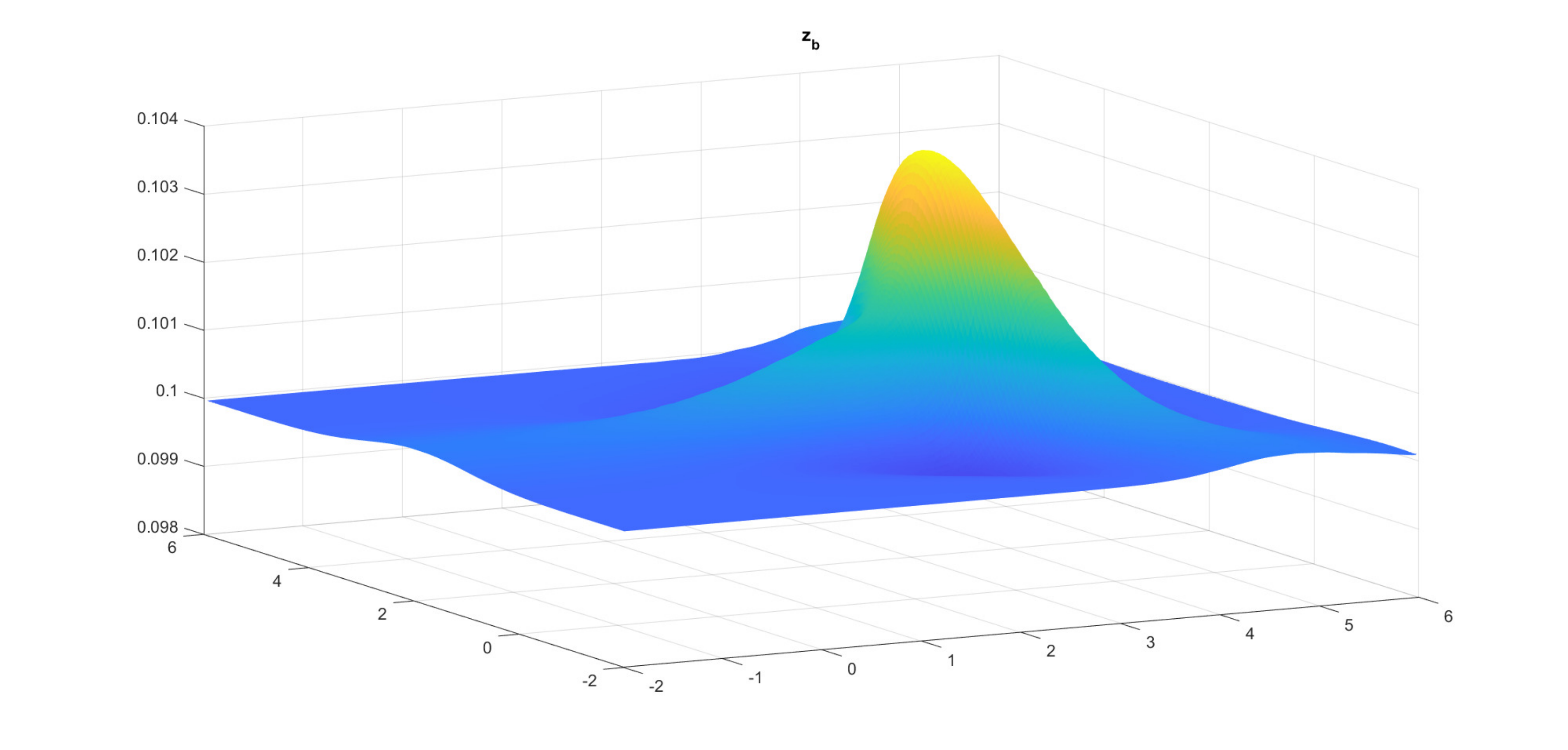}
    	\vspace{-0.8 cm}
    	\caption{Test \ref{Test_2_b}: (2D Exner conical sediment). Numerical solution for sediment layer  for the 2D Exner model on the square $[-2,6]\times[-2,6]$ using a $300\times300$ mesh points at time $t=2500$ with CFL$ = 12.$  }
    	\label{Test_2_2_con}
    \end{figure}
    \begin{figure}[!ht]
    	\centering
    	\includegraphics[height=0.36\textheight]{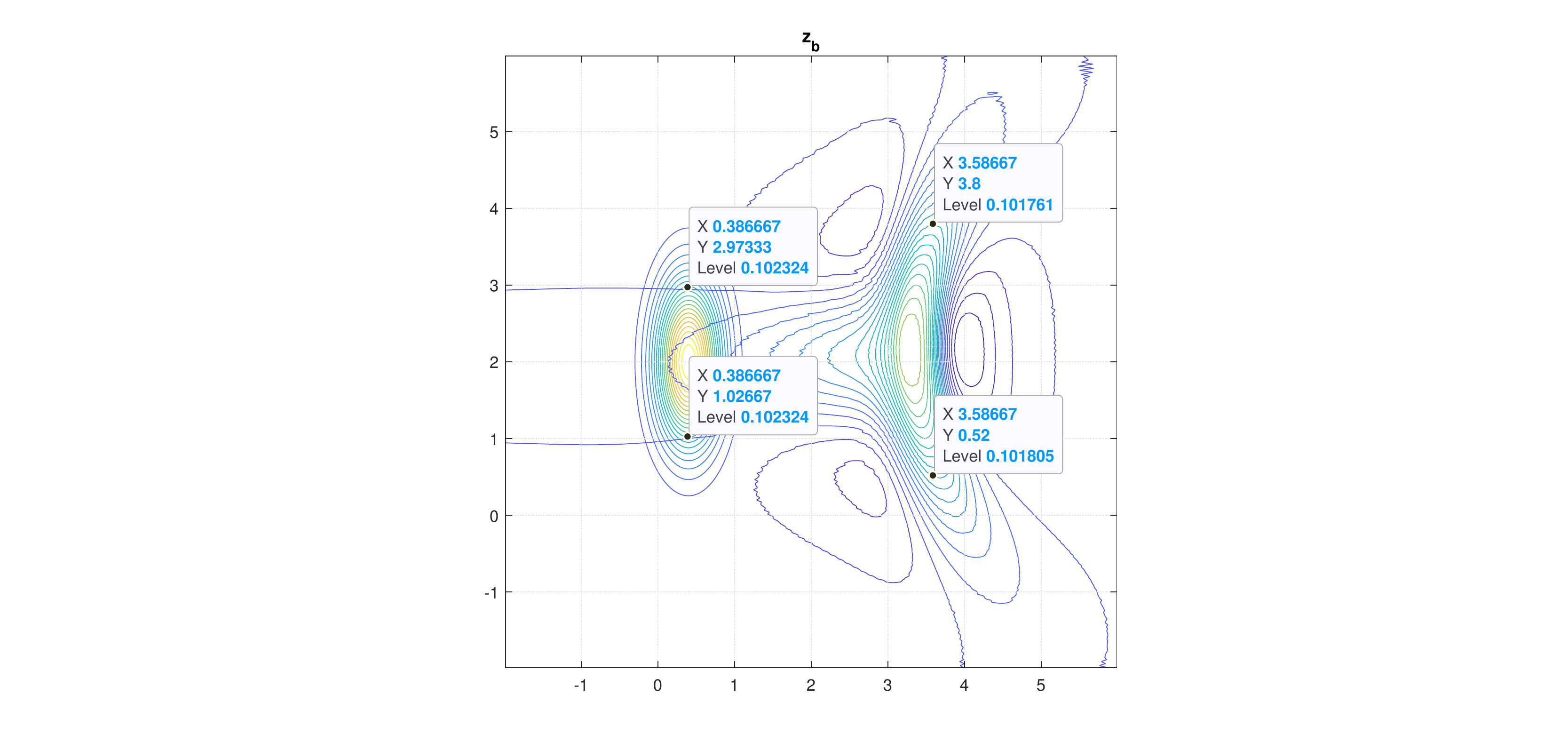}
    	\vspace{-0.8 cm}
     	\caption{Test \ref{Test_2_b}: (2D Exner conical sediment). Contour plot of initial and finale solution for sediment layer of 2D Exner model on the square $[-2,6]\times[-2,6]$ using a $300\times300$ mesh points at time $t=2500$ with CFL$ = 12.$ The lighted point are used to compute the spread angle.  }
    	\label{Test_2_2_cont}
    \end{figure}
     The numerical results are obtained with the second order semi-implicit scheme introduced on Section \ref{2D_Sec_ord_scheme} on the square $[-2,6]\times[-2,6]$ adopting a $300\times300$ mesh points, CFL$=12$ at time $t = 2500.$ Figure \ref{Test_2_2_con} shows the numerical evolution of the sediment layer. In particular Figure \ref{Test_2_2_cont} show initial and final solution and the relative points that we used to compute the spread angle. In this case, we consider a contour with 20 level and we compute the angle looking at level 8 from the bottom. The spread angle $\gamma \approx 23.14\degree .$ 
     
    \begin{figure}[!ht]
    	\centering	
    	\hspace{-1.2cm}
    	\includegraphics[scale = 0.16]{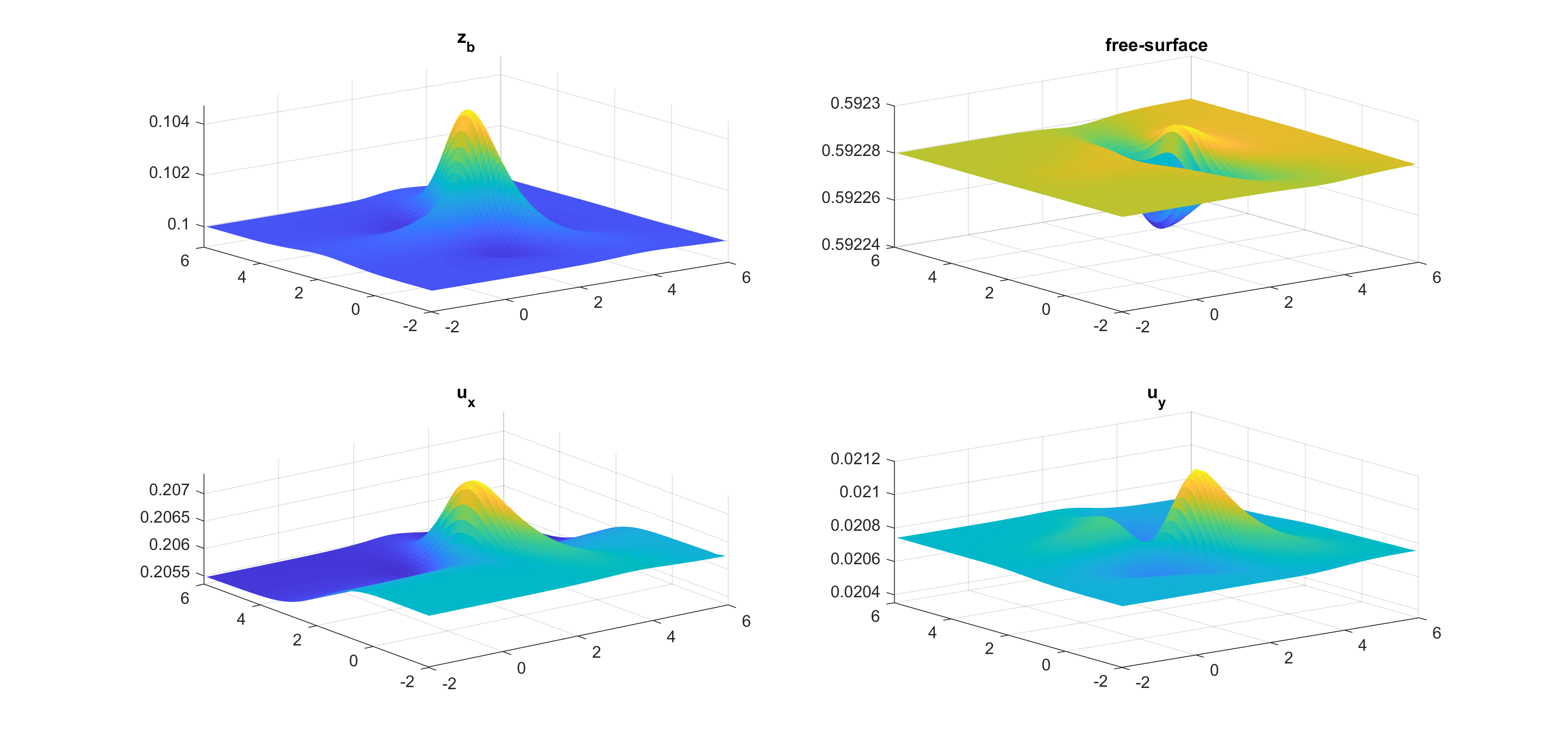}
    	\vspace{-0.8 cm}
    	\caption{Test \ref{Test_2_b}: (2D Exner conical sediment). Numerical solutions for sediment layer (up-left); free surface (up-right); velocity on $x-$direction (down-left) and velocity on $y-$direction (down-right) for the 2D Exner model on the square $[-2,6]\times[-2,6]$ using a $100\times100$ mesh points at time $t=300$ with CFL$ = 6.$  }
    	\label{Test_3_2_con}
    \end{figure}
    As last experiment we consider a fully two-dimensional case where initial conditions are so set: $\eta_0(x,y,0) = 0.6,$ $b(x,y) = 0,$ $m(x,y,0) = 0.2,$ $n(x,y,0)=0.02$ and 
    \begin{equation}
        \label{zb0_1_1} 
        z_{b_0}(x,y) = 0.1 + 0.006e^{-\frac{(x-0.4)^2}{0.4^2} - (y-3)^2}
    \end{equation}
    The numerical results are obtained with the second order semi-implicit scheme introduced on Section \ref{2D_Sec_ord_scheme} on the square $[-2,6]\times[-2,6]$ adopting a $100\times100$ mesh points, CFL$=6$ at time $t = 300.$ Figure \ref{Test_3_2_con} shows the numerical solutions of the sediment layer (up-left); free surface (up-right); velocity on $x-$direction (down-left) and velocity on $y-$direction (down-right). We can see that the flow on $y-$direction have and active role since a little phenomena appear on the velocity of $x-$direction. This behaviour is not visible in the sediment layer or free surface.
    
    The second order semi-implicit scheme is able to accurately compute  the evolution of the sediment.

\section{Conclusion}
The aim of this work is the development of semi-implicit schemes for the 1D and 2D Shallow-water Exner model. The objective was to drastically improve the  efficiency in the computation of the evolution of the sediment by treating water waves implicitly, thus allowing much larger time steps than the one permitted by standard CFL condition on explicit schemes, under the hypothesis of small Froude number and weak interaction between the fluid and the sediment layer. 

In particular, one of the objectives of this work is to check that even if we do not resolve the small time scale of the waves, still the semi-implicit method is able to correctly capture the sediment evolution. To this purpose, a simplified  scalar model, in which the flow is quasi-stationary, has been considered. As expected, there exists a very good agreement between the solution of the scalar equation and the full system, under the weak-interaction assumption. Furthermore, the long-time behaviour of the sediment has been checked even in presence of under-resolved fast water waves and the effects of these on the sediment have been analysed numerically.

There are still a few things that require improvements and generalizations:
\begin{itemize}
    \item In the current method the stability condition on the time step is determined by the fluid velocity, which may be much larger than the sediment wave speed. We plan to construct a scheme in which the CFL condition is based on the sediment wave rather than on the fluid velocity.
    \item More realistic Exner models can be considered. Other equations for the evolution of the sediment can be considered making the simulations more sophisticated and obtaining numerical results that are more consistent with the experimental data. 
    \item Different IMEX schemes should be explored and studied to improve accuracy and efficiency of the approach.
\end{itemize}

\section*{Acknowledgements} 
 This research has received funding from the European Union’s Horizon 2020 research and innovation program, under the Marie Sklodowska-Curie grant agreement No 642768 and from the European Union’s NextGenerationUE – Project: Centro Nazionale HPC, Big Data e Quantum Computing, “Spoke 1” (No. CUP E63C22001000006). E. Macca was partially supported by GNCS No. CUP E53C22001930001 Research Project “Metodi numerici per problemi differenziali multiscala: schemi di alto ordine, ottimizzazione, controllo”. E. Macca and G.Russo would like to thank the Italian Ministry of Instruction, University and Research (MIUR) to support this research with funds coming from PRIN Project 2017 (No. 2017KKJP4X entitled “Innovative numerical methods for evolutionary partial differential equations and applications”). M.J Castro research has been partially supported by the Spanish Government and FEDER through the coordinated Research project RTI2018-096064-B-C1, the Junta de Andaluc\'ia research project P18-RT-3163, the Junta de Andalucia-FEDER-University of M\'alaga research project UMA18-FEDERJA-16 and the University of M\'alaga. E. Macca and G. Russo are members of the INdAM Research group GNCS.

\bibliographystyle{acm}
\bibliography{biblio}

\end{document}